\documentclass[a4paper,11pt]{article}

\usepackage{graphicx}
\usepackage{latexsym}
\usepackage{amsmath}
\usepackage{amsthm}
\usepackage{amsfonts}
\usepackage{amssymb}
\usepackage{enumerate}

\usepackage[authoryear]{natbib}
  4  \oddsidemargin-20pt \evensidemargin-20pt
\allowdisplaybreaks[4]

\newcommand{\mockalph}[1]{}

\newcommand{\lfrf}[1]{\lfloor {#1}\rfloor}
\newcommand\N{\mathbb{N}}

\newcommand\R{\mathbb{R}}
\newcommand\Z{\mathbb{Z}}
\newcommand\eps{\varepsilon}
\newcommand\Prob{\mathbb{P}}    
\newcommand\E{\mathbb{E}}     

\newcommand\Fc{\mathcal{F}}
\newcommand\Gc{\mathcal{G}}

\newcommand\Lc{\mathcal{L}}
\newcommand\Bb{\mathbb{B}}
\newcommand\Cb{\mathbb{C}}

\newcommand\Gb{\mathbb{G}}

\newcommand\Pb{\mathbb{P}}
\newcommand\Vb{\mathbb{V}}

\newcommand\Yb{\mathbb{Y}}
\newcommand\Kb{\mathbb{K}}

\newcommand\Df{\mathbf{D}}
\newcommand\Rf{\mathbf{R}}
\newcommand\Cf{\mathbf{C}}
\newcommand\Uf{\mathbf{U}}
\newcommand\Vf{\mathbf{V}}

\DeclareMathOperator{\Cov}{Cov}

\DeclareMathOperator{\BL}{BL}

\def\cov{{\mbox{cov}}}

\newcommand\weak{\ \rightsquigarrow\ }

\newcommand\weakPi[1]{ \overset{\Prob}{\underset{ #1}{\,\rightsquigarrow}\,} }

\newcommand\Pconv{ \stackrel{ \Prob}{\rightarrow} }

\newtheoremstyle{normal}
{2ex}               
{3ex}               
{}                  
{}                  
{\bfseries} 
{}                  
{2pt}   
{\thmname{#1}\thmnumber{ #2.} \thmnote{(#3)}}

\newtheoremstyle{italic}
{2ex}
{3ex}
{\itshape}
{}
{\bfseries} 
{}
{2pt}
{\thmname{#1}\thmnumber{ #2.} \thmnote{(#3)}}

\theoremstyle{normal}
\newtheorem{definition}{Definition}[section]
\newtheorem{remark}[definition]{Remark}
\newtheorem{example}[definition]{Example}

\theoremstyle{italic}
\newtheorem{theorem}[definition]{Theorem}

\newtheorem{lemma}[definition]{Lemma}
\newtheorem{prop}[definition]{Proposition}

\begin{document}

\author{Stanislav Volgushev and Xiaofeng Shao\thanks{Stanislav Volgushev is Postdoctoral Researcher, Department of Mathematics, Ruhr University Bochum, 44780 Bochum, Germany. Volgushev's research was supported by the DFG grant Vo1799/1-1. Xiaofeng Shao is Associate Professor, Department of Statistics, University
of Illinois at Urbana-Champaign, Champaign, IL, 61820, USA. Shao's research is supported in part by NSF grant DMS-11-04545.  This research was conducted while Volgushev was visiting the University of Illinois at Urbana-Champaign. He would like to thank the people at the ´Statistics and Economics departments for their hospitality. The authors would also like to thank Axel B{\"u}cher and Holger Dette for comments on a preliminary version of this manucript.
Emails: stanislav.volgushev@rub.de, xshao@illinois.edu.}}
\title{A general approach to the joint asymptotic analysis of statistics from sub-samples}
\maketitle

\begin{abstract}
In time series analysis, statistics based on collections of estimators computed from subsamples play a crucial role in an increasing variety of important applications. Proving results about the joint asymptotic distribution of such statistics is challenging since it typically involves a nontrivial verification of technical conditions and tedious case-by-case asymptotic analysis. In this paper, we provide a novel technique that allows to circumvent those problems in a general setting. Our approach consists of two major steps: a probabilistic part which is mainly concerned with weak convergence of sequential empirical processes, and an analytic part providing general ways to extend this weak convergence to functionals of the sequential empirical process. Our theory provides a unified treatment of asymptotic distributions for a large class of statistics, including  recently proposed self-normalized statistics and  sub-sampling based p-values. In addition, we comment on the consistency of bootstrap procedures and obtain general results on compact differentiability of certain mappings that seem to be of independent interest.
\end{abstract}



\section{Introduction and Motivation}

In time series analysis, a large class of statistics can be expressed as  smooth functions of
estimators computed on  consecutive portions (i.e., subsamples) of data. Since time series observations are naturally ordered by time, the use of such statistics has been a common theme in time series inference and examples are abundant in areas such as sequential monitoring [\cite{CSW1995}, \cite{AHR2009}], retrospective change point detection [\cite{CH1997}, \cite{P2006}] and subsampling-based inference [\cite{PR1994}, \cite{PRW1999}], among others.  More recent examples include the self-normalized (SN, hereafter) statistics [\cite{shao2010}],  a new SN-based test statistic for change point detection [\cite{shzh2010}]  and the p-value of the subsampling-based inference under the fixed-b asymptotics [\cite{shpo2013}]. To obtain the asymptotic distributions of statistics of such kind, a traditional approach is to express the estimator as a sum of three parts, including the parameter, an average of influence functions, and a remainder term, followed by certain assumptions that ensure  asymptotic negligibility of remainder terms and a routine analysis of the leading term which is of linear form. For many statistics of practical interest, theoretical analysis based on this approach can be quite challenging and tedious. In particular verifying the negligibility of remainder terms can be technically involved, since it requires a careful case-by-case study. The situation is further complicated by the fact that in time series settings, the underlying data are dependent. The aim of the present paper is to provide a general approach which allows to easily obtain the asymptotic distribution of statistics based upon infinite collections subsample estimates without long and tedious arguments. 

In statistical applications, many important statistics can be expressed as smooth [more precisely: compactly differentiable] functionals of simple quantities such as the empirical distribution function. The analysis of the asymptotic properties of such statistics in the non-sequential setting can be elegantly performed in two distinct steps: an analytic part which consists in establishing the smoothness of the functional and a probabilistic part that is concerned with the analysis of the underlying quantity. One of the many appealing features of such an approach lies in the fact that the analytic properties need to be established only once. Moreover, quantities such as the empirical distribution function are often rather well analyzed for a wide range of data types. This approach has been successfully applied to the analysis of quantiles [\cite{dogi1992}], survival data [\cite{gijo1990}], copulas and scalar measures of dependence [\cite{ferradweg2004,buvo2011}] and to the setting of dependent data.

A slightly more formal description of the situation above is as follows. Assume that we have a collection of estimators, say $(\hat x_{n,\kappa})_{\kappa\in K}$ of a quantity $x$. A classical example of such a collection is given by estimators computed from various fractions of the sample $X_1,...,X_n$. For illustration purposes, assume that $x$ is the distribution function, $\hat x_{n,\kappa}$ denotes the empirical distribution function computed from $X_1,...,X_{\lfrf{n\kappa}+1}$ and $K = [0,1]$. Also, assume that the parameter of interest, say $\theta$, can be expressed as $\phi(x)$ where $\phi$ denotes some functional. For example, it is possible to express the copula as a functional of the cumulative distribution function. If the map $\phi$ is compactly differentiable, the asymptotic distribution of a suitably normalized version of $\phi(\hat x_{n,\kappa})$ for fixed $\kappa$ can be derived from a corresponding result for $\hat x_{n,\kappa}$. More precisely, denoting by $\alpha_n$ a sequence diverging to infinity and by $w(\kappa)$ a weight function, weak convergence of $\alpha_nw(\kappa)(\hat x_{n,\kappa}-x)$ in a suitable function space implies weak convergence of $ \alpha_nw(\kappa)(\phi(\hat x_{n,\kappa}) - \phi(x))$ \textit{for a finite collection of fixed values of $\kappa$}. However, in many important applications the \textit{joint} weak convergence of the whole \textit{collection} $\Vb_n := \Big(\alpha_nw(\kappa)(\phi(\hat x_{n,\kappa}) - \phi(x)\Big)_{\kappa\in K}$ in a suitable functional sense is required. For the purpose of illustration, consider the following simple example. 

\begin{example}
For the sake of concreteness, assume that we observe data, say $(Y_i = T_i\wedge C_i,\delta_i=I\{Y_i=T_i\})_{i=1,...,n}$ from a censored time series [here, $T_i$ denote survival times, $C_i$ censoring times and $\delta_i$ denote censoring indicators] and want to test if there is a change in the location parameter of the marginal distribution $F_i$ of $T_i$. A general way to quantify the location of censored observations, that is well-defined even under heavy censoring, is provided by the median. Typical test statistics for the null hypothesis of a constant median are based on comparing the medians of the Kaplan-Meier estimators which are computed from portions of the data. For simplicity, assume that the estimator $\hat m_\kappa$ with $\kappa \in [0,1]$ is based on the data $(Y_i,\delta_i)_{i=1,...,\lfrf{n\kappa}\vee 1}$. A simple test statistic for the null hypothesis of a constant median is given by $\sup_{\kappa\in[0,1]} w(\kappa)|\hat m_\kappa - \hat m_1|$ with $w$ denoting a suitable weighting function. In order to derive the null distribution of our test statistic, we would typically establish a process convergence result for $\sqrt{n}w(\kappa)(\hat m_\kappa - \hat m_1)$ viewed as element in the space $D[0,1]$ and apply the continuous mapping theorem.
Classical results on compact differentiability [see Example \ref{ex:quant} and \ref{ex:km}] imply that the median of the Kaplan-Meier estimator can be represented as a compactly differentiable functional of the two empirical (sub-)distribution functions $\hat H_{0,\lfrf{n\kappa}}(y) := \lfrf{n\kappa}^{-1}\sum_{i=1}^{\lfrf{n\kappa}} \delta_i I\{Y_i \leq y\}$ and $\hat F_{Y,\lfrf{n\kappa}}(y) := \lfrf{n\kappa}^{-1}\sum_{i=1}^{\lfrf{n\kappa}} I\{Y_i \leq y\}$. If we want to apply the classical delta-method to derive the process asymptotics of $\sqrt{n}w(\kappa)(\hat m_\kappa - \hat m_1)$ we are faced with two problems: first, we need process convergence of suitably normalized versions of $\hat H_{0,\lfrf{n\kappa}}(y) - \E[H_{0,\lfrf{n\kappa}}(y)]$ and $\hat F_{Y,\lfrf{n\kappa}}(y) - \E[F_{Y,\lfrf{n\kappa}}(y)]$. Second, as we shall argue below, the classical delta method does \textit{not} provide results on weak convergence of the quantity $\sqrt{n}w(\kappa)(\hat m_\kappa - \hat m_1)$ as a \textit{process} indexed in $\kappa$.
\end{example}

Returning to a more general setting, we can say that the classical delta method and a large collection of results on the behavior of general empirical processes allow to establish weak convergence results for a wide class of statistics as long as we consider a fixed, finite collection of values $\kappa$. Informally, we call this the 'non-sequential' case. However, the tools available to date do not allow the same conclusion when we are interested in collections of sub-samples, or, stated informally, in the 'sequential' case. The fundamental aim of the present article is thus to provide general ways of importing the tools mentioned above from the 'non-sequential' into the 'sequential' setting.  

For example, let us consider what we would need to apply a delta method in the 'sequential' case if only compact differentiability of the map $\phi$ in the 'non-sequential' case is available. Essentially, such an approach would require us to show compact differentiability of the map
\[
\Phi: (h_\kappa)_{\kappa\in K} \mapsto \Big(w(\kappa)\phi\Big(\frac{h_\kappa}{w(\kappa)}\Big)\Big)_{\kappa\in K}
\]
viewed as a map between suitable metric spaces since we can write
\[
\Vb_n = \alpha_n\Big(\Phi\Big( (w(\kappa)\hat x_{n,\kappa})_{\kappa\in K}\Big) - \Phi\Big( (w(\kappa)\hat x_{n,\kappa})_{\kappa\in K}\Big)\Big).
\]
Given the fact that a large amount of important maps $\phi$ that are known to be compactly differentiable, we would like to make use of this information in the sequential setting. A natural question to ask thus is: given compact differentiability of $\phi$, what can we say about compact differentiability of $\Phi$? As we shall see in Section \ref{sec:genana}, such an implication does not hold in full generality, see in particular Example \ref{ex:counter} and the discussion preceding it. At the same time, we obtain a positive result if we additionally assume that the map $\phi$ possesses certain boundedness properties. Additionally, even when compact differentiability of $\Phi$ fails, there still are many relevant settings where additional arguments can be applied to obtain the desired weak convergence of $\Vb_n$. In fact, in Section \ref{sec:genana} we show that, given weak convergence of $\Yb_n := \Big(\alpha_nw(\kappa)(\hat x_{n,\kappa} - x)\Big)_{\kappa\in K}$, we can derive properties of $\Vb_n$ in a very general setup. Additionally, some general results on compact differentiability that seem to be of independent interest can be found in Section \ref{sec:genhd}.\\
Another fundamental question that needs to be taken care of before we can apply the functional delta method is the weak convergence of the process $\Yb_n$. In fact, results on weak convergence of $\Yb_n$ in settings where the data $X_1,...,X_n$ are allowed to be dependent are limited. A summary of available results as well as new insights providing considerable extensions of those findings are collected in Section \ref{sec:genpro}.\\
Finally, in Section \ref{sec:app}, we illustrate how the general results presented in Section \ref{sec:gen} can be applied to obtain new insights regarding the properties of recently proposed methods including  self-normalization and generalizations thereof [Section \ref{sec:appsn}], fixed-b corrections for sub-sampling methods [Section \ref{sec:appfb}], and SN-based testing procedures for change-points [Section \ref{sec:appcp}]. Some comments on the applicability of our results to bootstrap methods are also provided.

\section{General results}\label{sec:gen}
We begin by introducing some relevant notation.
For arbitrary sets $\Fc_1,...,\Fc_J,K$ define the vector space
\[
\Lc^\infty(\Fc_1,...,\Fc_J;K) := \Big\{(H_{1,t},...,H_{J,t})_{t\in K}\Big| H_{j,t}\in\ell^\infty(\Fc_j)\forall j,t \ \ \sup_t \sup_j \sup_{f\in \Fc_j} |H_{j,t}(f)| < \infty\Big\}
\]
with norm
\[
\|(H_{1,t},...,H_{J,t})_{t\in K}\|_\Lc := \sup_t \sup_j \sup_{f\in \Fc_j} |H_{j,t}(f)|.
\]
Note that $\Lc^\infty(\Fc_1,...,\Fc_J;K)$ can be identified with $\ell^\infty(K \times\Fc_1)\times...\times \ell^\infty(K \times\Fc_J)$ by considering the relation
\[
(H_{1,t},...,H_{J,t})_{t\in K}\in \Lc^\infty(\Fc_1,...,\Fc_J;K) \quad \leftrightarrow \quad \Big((t,f)\mapsto H_{1,t}(f),...,(t,f)\mapsto H_{J,t}(f) \Big).
\]
By the definition of $\Lc^\infty(\Fc_1,...,\Fc_J;K)$, we have $\sup_t\sup_f |H_{j,t}(f)|<\infty$ for all $j=1,...,J$ so that the maps $(t_j,f_j)\mapsto H_{j,t_j}(f_j)$ are indeed bounded and thus elements of $\ell^\infty(K \times \Fc_j)$. In particular, if the product space $\ell^\infty(K \times\Fc_1)\times...\times \ell^\infty(K \times\Fc_J)$ is equipped with the maximum norm $\|(x_1,...,x_J)\|_{\max} := \max_j\|x_j\|_\infty$ induced by the supremum norms on its components, the identification given above is an isometry, that is
\[
\|(H_{1,t},...,H_{J,t})_{t\in K}\|_\Lc = \Big\|\Big((t,f)\mapsto H_{1,t}(f),...,(t,f)\mapsto H_{J,t}(f) \Big)\Big\|_{\max}.
\]
Weak convergence in $\Lc^\infty(\Fc_1,...,\Fc_J;K)$ is henceforth understood as weak convergence in the Hoffmann-J{\o}rgensen sense in the space $\Lc^\infty(\Fc_1,...,\Fc_J;K)$ as a subspace of $\ell^\infty(K \times\Fc_1)\times...\times \ell^\infty(K \times\Fc_J)$ [see \cite{vandervaart1996}, Chapters 1.4 and 1.5 for more details].

\begin{remark}
In most situations, the sets $\Fc_1,...,\Fc_J$ can be viewed as subsets of $\R^d$. For example, the empirical distribution function $( n^{-1}\sum I\{X_i \leq y\})_{y\in \R^d}$ of a sample of $d$-dimensional random variables $X_1,...,X_n$ is naturally indexed by the set $\R^d$. Another approach that fits nicely into the empirical process setting and will play a central role in Section \ref{sec:genpro}, is to consider classes of functions $\{x\mapsto f(x)|f\in \Fc_j\}$. In this setting, the empirical process can be elegantly written as $\Big( n^{-1}\sum_{i=1}^n f(X_i) - \E f(X_i)\Big)_{f\in \Fc_j}$, see \cite{vandervaart1996} for examples. For example, the empirical distribution function can also be viewed as element of $\ell^\infty(\Fc)$ with $\Fc$ denoting the collection of indicators of rectangles, that is $\Fc = \{x \mapsto I\{x\leq y\}|y\in\R^d\}$. By identifying the function $x \mapsto I\{x\leq y\}$ with the point $y \in \R^d$ we obtain a way to index $\Fc$ by $\R^d$ and vice versa. In most of the following theoretical developments, the form of $\Fc_j$ will be arbitrary unless explicitly specified otherwise.
\end{remark}

As discussed previously, the asymptotic analysis of statistics based on the process $\Vb_n$ can be performed by considering two distinct questions: the stochastic properties of $\Yb_n$ and the analytic properties of the map $\phi$.
Both questions will be addressed in this section in a general setting. In section \ref{sec:genana}, we present our analytic considerations. An overview of existing results regarding the stochastic part as well as their extension will be considered in section \ref{sec:genpro}. Finally, some general results on compact differentiability that seem to be of independent interest are provided in section \ref{sec:genhd}.

\subsection{Analytic considerations}\label{sec:genana}
This section is primarily concerned with the following questions: given a collection of estimators $(\hat y_{n,s,t})_{(s,t) \in K}$ such that for fixed $(s,t)\in K$ each $\hat y_{n,s,t}$ is an element of $\ell^\infty(\Fc_1)\times...\times\ell^\infty(\Fc_J)$ with $K\subset \Delta := \{(s,t)\in[0,1]^2|s \leq t\}$, a smooth (in a suitable sense) map $\phi: \ell^\infty(\Fc_1)\times...\times\ell^\infty(\Fc_J) \to \ell^\infty(\Gc_1)\times...\times\ell^\infty(\Gc_L)$, and weak convergence of the process [$\alpha_n$ denotes some deterministic sequence diverging to infinity]
\[
\Yb_n(s,t,f_1,...,f_J) := (t-s)\alpha_n(\hat y_{n,s,t}(f_1,...,f_J) - x(f_1,...,f_J))
\]
viewed as element of $\Lc^\infty(\Fc_1,...,\Fc_J;K)$, what can we say about weak convergence of $\Vb_n$ where
\[
\Vb_n(s,t,g_1,...,g_L) := (t-s)\alpha_n(\phi(\hat y_{n,s,t})(g_1,...,g_L) - \phi(x)(g_1,...,g_L))
\]
as element of $\Lc^\infty(\Gc_1,...,\Gc_L;K)$? And what can we say about bootstrap validity for $\Vb_n$ given a valid bootstrap procedure for $\Yb_n$?\\

For instance, consider the situation where we have a sample $X_1,...,X_n$. Assume that the quantity $x$ can be represented as $x = \Big((\E[f(X)])_{f\in \Fc_1},....,(\E[f(X)])_{f\in \Fc_J}\Big)$ for some classes of functions $\Fc_1,...,\Fc_J$, see the examples below. A prime example for the quantity $\hat y_{n,s,t}$ is given by the estimator computed from the sub-sample $X_{\lfrf{ns}+1},...,X_{\lfrf{nt}}$, that is
\begin{equation} \label{eq:xnt}
\hat y_{n,s,t} := \Big(\Big(\frac{1}{\lfloor nt \rfloor-\lfrf{ns}}\sum_{i=\lfrf{ns}+1}^{\lfloor nt \rfloor} f(X_i) \Big)_{f\in\Fc_1},...,\Big(\frac{1}{\lfloor nt \rfloor-\lfrf{ns}}\sum_{i=\lfrf{ns}+1}^{\lfloor nt \rfloor} f(X_i) \Big)_{f\in\Fc_J}\Big)
\end{equation}
where the empty sum is defined as zero and we set $'0/0 = 0'$ to take care of the case $\lfrf{ns}=\lfrf{nt}$.\\

Regarding the smoothness of $\phi$, we impose the following condition
\begin{itemize}
\item[(C)] The map
\[
\phi: \ell^\infty(\Fc_1)\times...\times\ell^\infty(\Fc_J) \supset D_\phi \rightarrow R_\phi \subset \ell^\infty(\Gc_1)\times...\times\ell^\infty(\Gc_L).
\]
is compactly differentiable at $x$ tangentially to $V\subset \ell^\infty(\Fc_1)\times...\times\ell^\infty(\Fc_J)$. Additionally, $0 \in V$ as well as $f\in V \Rightarrow cf \in V$ for all $c>0$.
\end{itemize}

In the 'classical' setting, compact differentiability is known to provide a good balance between strength of the differentiability concept that is needed for establishing a general functional delta method and the number of statistically relevant functionals that can actually be shown to be compactly differentiable. See \cite{vandervaart1996}, Chapter 3.9 for a more detailed discussion of this topic. Two particular examples are discussed below. Of course, there exists a vast collection of further examples [copulas, dependence measures, M- and L-estimators to name just a few] that are equally important but not discussed here because of space considerations. For a more detailed list we refer the interested reader to Chapter 3.9 in \cite{vandervaart1996} and the recent paper by \cite{gazh2011}.

\begin{example}\label{ex:quant}
\textbf{Empirical quantiles}\\
Consider the class of functions $\Fc: \{y \mapsto I\{y\leq t\}|t \in \R\}$. In this case, $\hat y_{n,s,t}$ is simply the empirical distribution function of the sub-sample $X_{\lfrf{ns}+1},...,X_{\lfrf{nt}}$. Consider the quantile map $\phi: F \mapsto (F^{-1}(\tau))_{\tau\in S}$ for some $S \subset (0,1)$ which now corresponds to $\Gc_1$. Applying this map to $\hat y_{n,s,t}$ yields collections of empirical quantiles of the sub-samples $X_{\lfrf{ns}+1},...,X_{\lfrf{nt}}$. Compact differentiability of the quantile map can be established under appropriate conditions, see Lemma 3.9.23 in \cite{vandervaart1996}.
\end{example}

\begin{example}\label{ex:km}
\textbf{Kaplan-Meier estimator}\\
Assume that we have right-censored observations of the form $(Y_i,\delta_i)_{i=1,...,n}$. It is a well-known fact that the Kaplan-Meier estimator $\hat F_{KM}$ [\cite{kame1958}], viewed as a map into the set of distribution functions on $[0,V]$ for a suitable $V<\infty$, is a compactly differentiable functional of the two functions
\[
\hat F_1(t) := \frac{1}{n}\sum_i \delta_iI\{Y_i\leq t\},\quad \hat F_Y(t) := \frac{1}{n}\sum_i I\{Y_i\leq t\},
\]
see Chapter 3.9 in \cite{vandervaart1996}. This suggests to consider the classes of functions
\[
\Fc_1 := \Big\{(y,\delta) \mapsto \delta I\{y\leq t\} \Big| t \in \R \Big\},\quad \Fc_2 := \Big\{(y,\delta) \mapsto I\{y\leq t\}\Big| t \in \R \Big\}.
\]
Combining this with the quantile mapping [see Example \ref{ex:quant}] easily allows to consider quantiles of the Kaplan-Meier estimator.
\end{example}

Regarding the process $\Yb_n$, we need the following assumption,
\begin{itemize}
\item[(W)] Assume that
\[
\Yb_n \weak \Yb \quad \text{in} \quad  \Lc^\infty(\Fc_1,...,\Fc_J;\Delta)
\]
where
\[
\Yb(s,t,f_1,...,f_J) := (\Yb_{1}(s,t,f_1),...,\Yb_{J}(s,t,f_J))
\]
and $\Yb_j, j=1,...,J$ are centered, Borel measurable processes.
\end{itemize}

\begin{remark}
A detailed discussion of condition (W) for estimators $\hat y_{n,s,t}$ of the form \eqref{eq:xnt} is provided in the next section.
However, there are interesting examples that go beyond the framework described above. For example, the classical empirical copula process [see \cite{rueschendorf1976}] is of the form
\[
\Cb_n^\circ(s,u) = \frac{1}{\sqrt{n}}\sum_{i=1}^{\lfloor sn \rfloor} \left( I\{ X_i \leq F_n^-(u)  \} - C(u) \right)
\]
where $F_n^-(u) := (F_{n1}^-(u_1),...,F_{nd}^-(u_d))$ denotes the vector of the generalized inverses of the marginal empirical distribution functions $F_{nj}(y) = n^{-1}\sum_i I\{X_{ij} \leq y\}$ and $C$ is the copula of the distribution of $X$. Note that $F_n^-(u)$ depends on all the data regardless of the value of $s$. The process $\Cb_n^\circ(s,u)$ can be coerced into the general framework of this section by considering the collection of estimators $\frac{1}{\lfrf{ns}}\sum_{i=1}^{\lfrf{ns}} I\{ X_i \leq F_n^-(u) \}$ indexed by $\Fc_1 := [0,1]^d$.
Weak convergence of the process $\Cb_n^\circ$ under weak assumptions on the copula with possibly dependent data was recently established by \cite{buvo2011}.
\end{remark}

The limit $\Yb$ in assumption (W) needs to satisfy certain technical conditions that are not very restrictive as we shall demonstrate later.

\begin{enumerate}
\item[(A1)] Assume that $\sup_{|s-s'| + |t-t'| \leq\delta}\sup_{j}\sup_{f_j\in \Fc_j}|\Yb_j(s,t,f_j) - \Yb_j(s',t',f_j)| = o_P(1)$ as $\delta\to 0$.
\item[(A2)] Define the set
\[
U_K :=\Big\{ (h_{s,t})_{(s,t)\in K}: h_{s,t} \in V\  \forall~(s,t)\in K, \sup_{(s,t)\in K}\|h_{s,t}\|<\infty \Big\}.
\]
Assume that the sample paths of $\Yb$ are in $U_K$ with probability one.
\end{enumerate}

Condition (A2) is non-restrictive in the sense that it is needed to apply the functional delta method to $\Vb_n(s,t,\cdot)$ for each fixed $(s,t)$. Assumption (A1) is needed for the application of the general compact differentiability result in Section \ref{sec:genhd}. As we shall discuss in the next section [see Remark \ref{rem:a1weak}], assumption (A1) is typically satisfied in a wide variety of practically relevant settings. Assumptions (W), (A1), (A2) are already sufficient to derive weak convergence of $\Vb_n$ if the set $K$ satisfies $\inf_{(s,t)\in K} |t-s| > 0$. Without this condition, (W), (A1), (A2) are not sufficient as the following example demonstrates.

\begin{example}\label{ex:counter}
Consider the map $\phi$ that takes a distribution function to its median and let $K = \{0\}\times[0,1] \subset \Delta$.
Define $\Fc_1 := \{x\mapsto I\{x\leq y\} | y\in \R\}$ and identify the functions $x\mapsto I\{x\leq y\} \in \Fc_1$ with $y\in \R$. Assume that $\hat y_{n,0,t}(y) = \frac{1}{\lfrf{nt}}\sum_{i=1}^{\lfrf{nt}}I\{X_i\leq y\}$. Consider a triangular scheme of data that is of the form $X_{jn} = n$ for $1 \leq j < n^{1/3}$ and $X_{jn}\sim U[0,1]$ i.i.d. for $n^{1/3}\leq j \leq n$. Elementary calculations show that $\Yb_n(0,\cdot,\cdot)$ converges weakly to the Kiefer-M{\"u}ller process $\Kb$ with covariance $\Cov(\Kb(t,y),\Kb(t',y')) = \min(t,t')(\min(y,y')-yy')$. On the other hand, setting $t = n^{-3/4}$ we have almost surely

\[
\Vb_n((0,n^{-3/4})) = n^{1/2}t(\phi(\hat y_{n,0,t}) - \phi(x)) = n^{-1/4}(n-1/2) \to \infty,
\]
and thus weak convergence of $\Vb_n$ can not hold.
\end{example}

The underlying problem in the above example is that due to the weighting with $t-s$, weak convergence of $\alpha_n(t-s)(\hat y_{n,t,s} -x)$ is not informative about $\hat y_{n,t,s}$ for values of $t-s$ that can be arbitrarily close to zero. Additional assumptions are needed to exclude this kind of behavior presented in the above example. It turns out that for this purpose the following condition is sufficient. As we shall discuss later, there are many situations where it is easily satisfied.

\begin{enumerate}
\item[(A3)] For any $k_n \to 0$ we have $\sup_{(s,t)\in K,|t-s|\leq k_n} (t-s)\|\phi(\hat y_{n,s,t})\| = o_P^*(1)$ where the asterisk denotes outer probability.
\end{enumerate}

\begin{remark}
Note that condition (A3) is automatically satisfied if $\sup_{(s,t)\in K}\|\phi(\hat y_{n,s,t})\| = O_P^*(1)$. This is trivially true for uniformly bounded maps $\phi$, which includes many interesting examples such as copulas, dependence measures or the Kaplan-Meier estimator (which per definition is a distribution function). Moreover for specific sets $K$, further conditions implying (A3) can be derived. See Remark \ref{rem:a3} in Section \ref{sec:genpro} for further details.
\end{remark}
Now we are ready to state the first main result of this section.

\begin{theorem} \label{theorem:main1}
For any compact $K\subset\Delta$, with $\inf_{(s,t)\in K}|t-s|\geq a > 0$ conditions (C), (W), (A1) and (A2) imply $\Vb_n\weak \Vb$ in $\Lc^\infty(\Gc_1,...,\Gc_L;K)$ where
\[
\Vb((s,t),g_1,...,g_L) := \Big(\phi_x' \Yb(s,t,\cdot)\Big)(g_1,...,g_L).
\]
If additionally (A3) holds, the assumption $\inf_{(s,t)\in K}|t-s|\geq a > 0$ can be dropped.
\end{theorem}

\begin{remark}\label{rem:clip}
Although assumption (A3) often holds, there are situations where verifying it can be very tedious or requires additional assumptions on the underlying data structure. For example, consider the setting where $K = \Delta$ and $\phi$ denotes the map that takes a distribution function to its median. In that case, assumption (A3) would require that $\frac{1}{n}\max_{i=1,...,n} |X_i| = o_P(1)$ since the median of one observation is the observation itself. Effectively, this places moment assumptions on $X$ that are not needed for the median from large samples to be well-behaved. A closer look at the proofs reveals that for any $\gamma \in (0,1)$ the following modified version of the process $\Vb_n$
\[
\tilde \Vb_n := \Big((t-s)I\{t-s \geq \alpha_n^{-\gamma}\}\alpha_n(\phi(\hat y_{n,s,t}) - \phi(x))\Big)_{(s,t)\in K}
\]
converges to the same limit $\Vb$ without assumption (A3) or the condition $\inf_{(s,t)\in K}|t-s|\geq a > 0$. In the applications discussed in Section \ref{sec:app}, the modification above essentially amounts to not using information from extremely small sub-samples. As the discussion above indicates, for certain sets $K$ this can be viewed as a robustification.
\end{remark}

\begin{remark}\label{rem:gencom}
A closer look at the proof of the above result shows that the special structure of $\Delta$ does not play a crucial role. In fact, the same approach yields a more general result. Let $(K,d_K)$ denote a general compact metric space. Assume that
\[
\Yb_n(\kappa,f_1,...,f_J) = w(\kappa)\alpha_n\Big(\hat y_{n,\kappa}(f_1,...,f_J) - x(f_1,...,f_J)\Big)
\]
and that $\Yb_n \weak \Yb \in \Lc(\Fc_1,...,\Fc_J;K)$ with a centered process $\Yb$ with $w(\cdot)$ denoting a bounded weight function. If additionally $\sup_{d_K(\kappa,\kappa') \leq\delta}\sup_{j}\sup_{f_j\in \Fc_j}|\Yb_j(\kappa,f_j) - \Yb_j(\kappa',f_j)| = o_P(1)$ and if the sample paths of $\Yb$, are, with probability one in the set 
\[
\tilde U_K := \Big\{ (h_{\kappa})_{\kappa\in K}: h_{\kappa} \in V\  \forall~\kappa\in K, \sup_{\kappa\in K}\|h_{\kappa}\|<\infty \Big\}
\]
it follows that with $\gamma \in (0,1)$ arbitrary
\[
\Big(w(\kappa)I\{w(\kappa) \geq \alpha_n^{-\gamma}\}\alpha_n(\phi(\hat y_{n,\kappa}) - \phi(x))\Big)_{\kappa\in K}  \weak \Big(\phi_x' \Yb(\kappa,\cdot)\Big)_{\kappa\in K}\quad \text{in} \quad \Lc^\infty(\Gc_1,...,\Gc_L;K)
\]
as long as $\inf_{\kappa\in K}|w(\kappa)|>0$. If additionally a modified version of condition (A3) holds, i.e. if for any $k_n \to 0$ we have $\sup_{\kappa\in K,|w(\kappa)|\leq k_n} w(\kappa)\|\phi(\hat y_{n,\kappa})\| = o_P^*(1)$, the weak convergence above holds without the assumption $\inf_{\kappa\in K}|w(\kappa)|>0$.
\end{remark}

Next, we discuss bootstrap procedures.
 In particular, consider the following bootstrap version of the quantity $\hat y_{n,s,t}$ defined in \eqref{eq:xnt}
\begin{equation}\label{eq:ynb1}
\hat y_{n,s,t}^b(s,t,f_1,...,f_J) := \Big(\frac{1}{\lfloor nt \rfloor - \lfloor ns \rfloor}\sum_{i=\lfloor ns \rfloor+1}^{\lfloor nt \rfloor} M_if(X_i),...,\frac{1}{\lfloor nt \rfloor - \lfloor ns \rfloor}\sum_{i=\lfloor ns \rfloor+1}^{\lfloor nt \rfloor} M_if(X_i) \Big) .
\end{equation}
with $M_1,...,M_n$ denoting random variables independent of the original sample $X_1,...,X_n$. The corresponding bootstrap version of the process $\Gb_n$ is given by
\begin{equation}\label{eq:ynb}
\Yb_n^b := (t-s)\alpha_n(\hat y_{n,s,t}^b - \hat y_{n,s,t})=: (\Yb_{n,1}^b,...,\Yb_{n,J}^b)
\end{equation}
Under suitable assumptions on the data and random variables $M_1,...,M_n$ a conditional version of assumption (W) holds. Specifically, assume that
\begin{itemize}
\item[(WB)] $\Yb_{n}^b$ weakly converges to $\Yb$ conditionally on the data in probability, or
\begin{equation*}
\label{as:condweak}
\Big( \Yb_{n,1}^b,...,\Yb_{n,J}^b\Big) \weakPi{M} \Big( \Yb_{1},...,\Yb_{J}\Big) \quad \mbox{in} \quad \Lc^\infty(\Fc_1,...,\Fc_J;K).
\end{equation*}
\end{itemize}
Here, weak convergence conditional on the data in probability ($\weakPi{M}$-convergence) is understood in the Hoffmann-J\o rgensen sense as defined in \cite{kosorok2008}, that is $\Yb_n^b\weakPi{M} \Yb$ if and only if
\begin{enumerate}
\item[(i)]  $\sup_{f\in \BL_1} \left| \E_M f(\Yb_n^b) - \E f(\Yb) \right| \rightarrow 0$ in outer probability,
\item[(ii)] $\E_M f(\Yb_n^b)^* - \E_M f(\Yb_n^b)_* \Pconv 0$ for all $f\in \BL_1$,
\end{enumerate}
where $\BL_1$ denotes the set of all Lipschitz-continuous functions $f:\Lc^\infty(\Fc_1,...,\Fc_J;K)\rightarrow\R$ that are uniformly bounded by $1$ and have Lipschitz constants bounded by $1$, and where the asterisks in (ii) denote measurable majorants (and minorants, respectively) with respect to the joint data $(X_1,\dots,X_n,M_1,\dots,M_n)$. Also, note that the map $(M_1,...,M_n) \mapsto \Yb_n^b$ is measurable conditionally on the original data $X_1,...,X_n$ outer almost surely [for fixed $X_1,...,X_n$, this mapping is Lipschitz-continuous] and thus we do not need to consider measurable majorants. Settings where results of this kind hold are discussed in the next section.\\
\\
The classical delta method for the bootstrap [see e.g. Theorem 12.1 in \cite{kosorok2008}] asserts that for a map $\phi$ that is compactly differentiable at $x$ with derivative $\phi'_x$ and additionally satisfies suitable measurability conditions, we have
\[
\alpha_n (t-s)(\phi(\hat y_{n,s,t}^b) - \phi(\hat y_{n,s,t})) \weakPi{M} \phi_x' \Gb \quad \mbox{in} \quad \ell^\infty(\Gc_1)\times...\times\ell^\infty(\Gc_L)
\]
for every fixed $(s,t)$. The next Theorem provides a generalization of this finding. More precisely, it states conditions that allow for a generalization of Theorem \ref{theorem:main1} to conditional weak convergence in $\Lc^\infty(\Gc_1,...,\Gc_L;K)$.

\begin{theorem} \label{th:boot1}
With the notation above, assume that (WB), (A1), (A2) and (C) hold. Then for any compact $K\subset \Delta$ with $\inf_{(s,t)\in K}|t-s|>0$ we have
\[
\Vb_n^b := \alpha_n (t-s)(\phi(\hat y_{n,s,t}^b) - \phi(\hat y_{n,s,t})) \weakPi{M} \phi_x' \Yb= \Vb  \quad \mbox{in} \quad \Lc^\infty(\Gc_1,...,\Gc_L;K).
\]
If additionally (A3) holds and $\sup_{(s,t)\in K,|t-s|\leq k_n} (t-s)\|\phi(\hat y_{n,s,t}^b)\| = o_P^*(1)$, the convergence holds for arbitrary compact $K\subset \Delta$.
\end{theorem}

\begin{remark}
Suitable modifications of the extensions discussed in Remark \ref{rem:clip} and Remark \ref{rem:gencom} continue to hold in the bootstrap setting. More precisely, we can replace sets $K\subset \Delta$ by arbitrary compact sets and the weighting $t-s$ with arbitrary bounded weighting functions $w$, in which case the assumption $\inf_{(s,t)\in K}|t-s|>0$ needs to be replaced by $\inf_{\kappa \in K}|w(\kappa)|>0$. Also, conditional weak convergence of $\tilde \Vb_n^b = (w(\kappa)I\{|w(\kappa)\}|>\alpha_n^{-\gamma}\}\alpha_n(\phi(\hat y_{n,s,t}^b) - \phi(\hat y_{n,s,t}))_{\kappa\in K}$ holds without assumption (A3) and the condition $\sup_{(s,t)\in K,|t-s|\leq k_n} (t-s)\|\phi(\hat y_{n,s,t}^b)\| = o_P^*(1)$ used in the above theorem.
\end{remark}

\subsection{Probabilistic considerations}\label{sec:genpro}

In this section, we focus our attention on the setting where $\Yb_n$ has a specific structure that typically arises in applications. More precisely, consider the multi-parameter sequential empirical process
\[
\Yb_{n} := \Big((t-s)\alpha_n(\hat y_{n,s,t} -x)\Big)_{(s,t)\in \Delta}
\]
where $\Delta := \{(s,t)\in [0,1]^2| s\leq t\}$ and the quantity
\[
\hat y_{n,s,t} := \Big(\Big(\frac{1}{\lfloor nt \rfloor - \lfloor ns \rfloor}\sum_{i=\lfloor ns \rfloor+1}^{\lfloor nt \rfloor} f(X_i) \Big)_{f\in\Fc_1},...,\Big(\frac{1}{\lfloor nt \rfloor - \lfloor ns \rfloor}\sum_{i=\lfloor ns \rfloor+1}^{\lfloor nt \rfloor} f(X_i) \Big)_{f\in\Fc_J}\Big)
\]
denotes an estimator for $x$ that is computed based on the sub-sample $X_{\lfloor ns \rfloor+1},...,X_{\lfloor nt \rfloor}$.
It turns out that conditions (W), (A1), (A2) in the previous section can be derived from simpler conditions that involve only a collection of 'classical' one-parameter sequential processes
\[
\Gb_n(t,f_1,...,f_J) := (\Gb_{n,1}(t,f_1),...,\Gb_{n,J}(t,f_J))
\]
where $\Gb_{n,j}(t,f) := t\alpha_n(\hat x_{n,t}^{(j)}(f)-x)$ and
\[
\hat x_{n,t} := \Big(\Big(\frac{1}{\lfloor nt \rfloor}\sum_{i=1}^{\lfloor nt \rfloor} f(X_i) \Big)_{f\in\Fc_1},...,\Big(\frac{1}{\lfloor nt \rfloor}\sum_{i=1}^{\lfloor nt \rfloor} f(X_i) \Big)_{f\in\Fc_J}\Big) =: (\hat x_{n,t}^{(1)},...,\hat x_{n,t}^{(J)}).
\]
Consider the assumptions
\begin{itemize}
\item[(W')] Assume that
\[
\Gb_n \weak \Gb \quad \text{in} \quad  \Lc^\infty(\Fc_1,...,\Fc_J;[0,1])
\]
where
\[
\Gb(t,f_1,...,f_J) := (\Gb_{1}(t,f_1),...,\Gb_{J}(t,f_J))
\]
and $\Gb_j, j=1,...,J$ are centered, Borel measurable processes.
\item[(A1')] Assume that $\sup_{|s-t|\leq\delta}\sup_{j}\sup_{f_j\in \Fc_j}|\Gb_j(t,f_j) - \Gb_j(s,f_j)| = o(1)$ as $\delta\to 0$.
\end{itemize}

The conditions above turn out to be sufficient for (W) and (A1).

\begin{prop}\label{prop:mult}
Under conditions (W') and (A1'), we have
\[
\Yb_n \weak \Yb \quad \text{in} \quad  \Lc^\infty(\Fc_1,...,\Fc_J;\Delta)
\]
where $\Delta = \{s,t \in [0,1]: s \leq t\}$ and
\[
\Yb(s,t,f_1,...,f_J) := \Gb(t,f_1,...,f_J) - \Gb(s,f_1,...,f_J).
\]
Moreover, $\Yb$ satisfies assumption (A1).
\end{prop}

\begin{remark}\label{rem:a1weak}
For many kinds of weakly dependent data [including, of course, the independent case], the process $\Gb$ is a vector of centered Gaussian processes with covariance of the form
\[
\Cov(\Gb(s,f_1,...,f_J),\Gb(t,g_1,...,g_J)) = (s \wedge t) K(f_1,...,f_J,g_1,...,g_J)
\]
for some uniformly bounded covariance kernel $K$. In this case, assumption (A1') holds.
To see this, note that under (A1') the process $\Gb$ has paths that are uniformly continuous with respect to the metric $\rho_2((t,f_1,...,f_J),(t',f_1',...,f_J')) := \E[(\Gb(t,f_1,...,f_J) - \Gb(t',f_1',...,f_J'))^2]$, see Example 1.5.10 in \cite{vandervaart1996}. The discussion at the beginning of Example 1.5.10 in \cite{vandervaart1996} thus yields the desired result. The special structure of $\Yb_n$ implies that its sample paths have the same property.
\end{remark}

\begin{remark}
There are interesting cases where condition (A1') holds for limiting processes that are non-Gaussian. More precisely, defining $\Fc_1 = [-\infty,\infty]$, the results in \cite{deta1989} imply weak convergence of the process $\Gb_n$ if the data $X_i$ exhibit long-range dependence. The limiting process, which can be non-Gaussian, is of the form $\Gb(t,y) = f(y)Z_m(t)$ with $f$ denoting a deterministic, uniformly bounded function and $Z_m$ a so-called m'th order Hermite-process [see \cite{deta1989} for a definition]. In particular, the sample paths of this process are H{\"o}lder-continuous [see \cite{matu2007}] and thus assumptions (A1') and (W') hold.
\end{remark}


\begin{remark}\label{rem:a3}
Consider the special case $K = \{0\}\times[0,1]$. In this case,
assumption (A3) is satisfied as soon as $\hat x_{n,t}$ is of the form given in (\ref{eq:xnt}) with the data $X_1,X_2,...$ stemming from a strictly stationary sequence and $\hat x_{n,1} \to x$ outer almost surely. To see this, note that under the assumptions discussed above we have $\sup_{t}\|\phi(\hat x_{n,t})\| = \max_{j=1,...,n} \|\phi(\hat x_{j,1})\|$ and that by Lemma \ref{lem:phi} together with the continuous mapping theorem $\|\phi(\hat x_{n,1})\| \to \|\phi(x)\|$ outer almost surely. This in turn implies that $(\sup_{n\geq 1}\|\phi(\hat x_{n,1})\|)^*$ [the asterisk denoting a measurable majorant] is bounded in probability. For results implying almost sure convergence in a very general setting, see \cite{adno2010} and the references cited therein.
\end{remark}

 For independent data, assumption (W') is known to hold as soon as the classes of functions $\Fc_1,...,\Fc_J$ are Donsker [see \cite{vandervaart1996}, Chapter 2.12.1]. For dependent data, much less is known. Available results are, to the best of our knowledge, limited to classes of functions of of the form $\Fc_1 = \{u\mapsto I\{u\leq y\}|y\in \R^d\}$ [the inequality is understood component-wise]. Here, results for $d>1$ are derived by \cite{sen1974} and \cite{ruschendorf1974} under $\phi-$mixing  and by \cite{yoshihara1975} and \cite{inoue2001} under strong mixing. \cite{beschho2009} considered the case $d=1$ under S-mixing, and derived a stronger result than weak convergence of the process. Finally, the paper by \cite{deta1989} contains a similar result for the class of functions $\Fc_1 = \{u\mapsto I\{u\leq y\}|y\in \R\}$ and long-range dependent data. To the best of our knowledge, nothing is known for general classes of functions.

Note that by Lemma 1.4.3 in \cite{vandervaart1996}, asymptotic tightness of $\Gb_n$ is equivalent to asymptotic tightness of $\Gb_{n,j}$ for all $j=1,...,J$. Thus, Problem 1.5.3 in the same reference implies that in order to obtain weak convergence of $\Gb_n$ to $\Gb$, we need to show that first $\Gb_{n,j}$ is asymptotically tight for all $j=1,...,J$ and second that the following condition holds
\begin{enumerate}
\item[(F)] For all finite collections $s_{i,j} \in [0,1], i=1,...,N, j=1,...,J$, $f_{ij} \in \Fc_j, i=1,...,N, j=1,...,J$ the collection $(\Gb_{j,n}(s_{ij},f_{ij}))_{j=1,..,J,i=1,...,N}$ converges weakly to $(\Gb_{j}(s_{ij},f_{ij}))_{j=1,..,J,i=1,...,N}$ in the usual $\R^{NJ}$-dimensional sense.
\end{enumerate}
There is a vast literature containing results that imply the finite-dimensional convergence (F), see \cite{demiso2002} and the references cited therein for an overview. Criteria establishing asymptotic tightness of the processes $\Gb_{n,j}$ for dependent data on the other hand are not as widely available, and one general result along those lines is provided below. This result is of independent interest. In particular, it can be used to verify condition (W') in a number of settings that have not been considered before.

\begin{theorem}
\label{th:donsker}
Assume that the process $\Gb_n$ is of the form $\Gb_n = t\alpha_n(\hat x_{n,t} - x)$ where $\hat x_{n,t}$ is defined in \eqref{eq:xnt} and the data $X_1,X_2,...$ come from a strictly stationary sequence. Assume that for each $j=1,...,J$ there exists a semi-metric $\rho_j$ on $\Fc_j$ which makes $\Fc_j$ totally bounded, and for each $j=1,...,J$ we have $\sup_{f\in \Fc_j} \E|f|^q < \infty$. Define ${\cal F}_{j,\delta} := \{f-g|f,g\in \Fc_j,\rho_j(f,g)\leq\delta\}$. Assume that the process $\Gb_{n,j}(1,\cdot)$ satisfies for some $q>2$ and $j=1,...,J$
\begin{eqnarray}
\label{eq:res1}
\lim_{\delta\downarrow 0}\limsup_{n\rightarrow\infty} \E^*\|\Gb_{n,j}(1,\cdot)\|_{{\cal F}_{j,\delta}}^q=0
\end{eqnarray}
[remember that the asterisk denotes outer expectation], that
\begin{eqnarray}
\label{eq:res3}
\max_{j=1,...,J}\sup_{n\in\N}\E^* \|\Gb_{n,j}(1,f)\|^q <\infty \quad \forall f\in \Fc.
\end{eqnarray}
and that for every $j$ the class of functions $\Fc_j$ has envelope $F_j$ which has finite $q$'th moment. Let condition (F) hold. Then $\Gb_n \weak \Gb$ in $\Lc^\infty(\Fc_1,...,\Fc_J;[0,1])$.
\end{theorem}

Condition (\ref{eq:res1})  has been established by \cite{anpo1994} for strongly mixing data, and inequality (3.1) in \cite{anpo1994} reveals that \eqref{eq:res3} holds under the same assumption. Moreover \cite{hagemann2012} established (\ref{eq:res1}) for stationary sequences with geometric moment contraction properties [see \cite{wush2004}], and the results in his appendix show that again \eqref{eq:res3} holds under the same assumptions.\\
\\
Next, consider bootstrap procedures. In the case of independent data, a mild assumption on the multipliers $M_i$ suffices. More precisely, assuming that $M_i$ are i.i.d., independent of the data $X_i$, and that $\int \sqrt{P(|M_1|>u)} du$ is finite [which follows if $M_1$ has finite moment of order $2+\eps$], the classes of functions $\Fc_1,...,\Fc_J$ being Donsker [see \cite{vandervaart1996}, page 81 for a definition of this property] implies (WB). To see this, note that by arguments similar to the ones given in the proof of Proposition \ref{prop:mult} it suffices to derive (WB) for the set $K = \{0\}\times[0,1]$. To do so, apply Lemma \ref{lem:bootappr} in the appendix where the approximating mappings $A_i$ and $A_{i,n}^b$ are defined through projections on piecewise constant functions, see the arguments in the proof of Theorem 1.5.6 in \cite{vandervaart1996}. Then assumption (i) of Lemma \ref{lem:bootappr} corresponds to conditional finite-dimensional convergence which can be established by arguments similar to those given in Lemma 2.9.5 in \cite{vandervaart1996}. Condition (ii) corresponds to tightness of the limit process $\Yb$. Condition (iii) follows from the unconditional asymptotic tightness of $\Yb_n^b$, which can be established by combining Theorem 2.12.1 and 2.9.2 in \cite{vandervaart1996}.\\
\\
Under dependence, much less is known about bootstrap validity for empirical processes, even in the non-sequential setting. For an overview of available results, see \cite{radulovic2009}.
In the sequential setting, some results along those lines were recently considered by \cite{buru2013} based on arguments from \cite{}. More precisely, those authors proposed to consider variables $M_{1,n},...,M_{n,n}$ from a triangular scheme that satisfy certain conditions [see assumptions A1-A3 in their paper]. In particular, the results in \cite{buru2013} imply (WB) for $K = \{0\}\times[0,1]$ under strong mixing conditions for the class of functions $\Fc = \{u\mapsto I\{u\leq w\}|w\in \R^d\}$. Moreover, using the techniques in that paper, in particular the Ottaviani type inequality [Lemma 1 in Appendix B of the corresponding paper], it should be possible to derive (WB) for $K = \{0\}\times[0,1]$ by combining arguments similar to those in the proof of Theorem 2.12.1 in \cite{vandervaart1996} with the Ottaviani-type inequality of \cite{buru2013} and results on the validity of bootstrap procedures in the non-sequential setting. For an overview of such results, see \cite{radulovic2009} and the references cited therein.

\subsection{A general result on (quasi) Hadamard differentiability}\label{sec:genhd}

This section contains an abstract result on compact differentiability that seems to be of independent interest. It plays a crucial role in the proofs of Theorems \ref{theorem:main1} and \ref{th:boot1}. The result in this section applies to both classical Hadamard differentiability [also known as compact differentiability], and the more general concept of quasi-Hadamard differentiability which was recently introduced by \cite{beza2010}. The main advantage of this more general approach is that it allows to apply a modified delta method in settings where the classical delta method fails, the simplest example being the mean. In particular, the distribution of U- and V-statistics and value-at-risk functionals can be derived in settings where the classical delta method fails. See \cite{beza2010,beza2012,bewuza2012} for further details. For the reader's convenience, we state the definition from \cite{beza2010}.

\begin{definition}[\cite{beza2010}]
Consider a metrized topological vector space $(R,d_R)$, a vector space $D$ with subsets $D_\phi,D_0\subset D, C_0\subset D_0$ and assume that $(D_0,d_D)$ is a metrized topological vector space. A map $\phi: D_\phi \to R$ is said to be \textit{quasi-Hadamard differentiable at $x\in D_\phi$ tangentially to $C_0\langle D_0\rangle$ with derivative $\phi_x'$} if for every $t_n \searrow 0$ and sequence $h_n \to h$ with $h_n\in D_0 \forall n$, $h\in C_0$ such that $x+t_nh_n \in D_\phi \forall n$ we have
\[
d_R(t_n^{-1}(\phi(x+t_nh_n) - \phi(x)),\phi_x'h) \to 0
\]
with $\phi_x'$ denoting a continuous, linear map $C_0 \to R$.
\end{definition}

Consider the following general setting.
\begin{itemize}
\item[(S)] Denote by $(R,d_R)$ a metrized topological vector space. Consider a second vector space $D$ with subsets $D_\phi,D_0 \subset D, C_0 \subset D_0$ and assume that $(D_0,d_D)$ is a metrized topological vector space. Let $\phi: D_\phi \to R$ be quasi Hadamard differentiable at $x$ tangentially to $C_0\langle D_0\rangle$ and denote the derivative by $\phi_x'$. Let $(K,d_K)$ be a compact metric space. Define the sets
\begin{eqnarray*}
\Df &:=& \Big\{(h_t)_{t\in K} \Big| h_t \in D~\forall t\Big\}
\\
\Rf &:=& \Big\{(h_t)_{t\in K} \Big| h_t \in R~ \forall t, \sup_{s,t} d_R(h_t,h_s) < \infty\Big\}
\\
\Df_\Phi &:=& \Big\{(h_t)_{t\in K} \Big| h_t \in D_\phi~ \forall t,\ \sup_{s,t} d_R(\phi(h_t),\phi(h_s)) < \infty\Big\}
\\
\Rf_\Phi &:=& \Big\{(h_t)_{t\in K} \Big| h_t \in D~ \forall t, \sup_{s,t} d_R(h_t,h_s) < \infty\Big\}
\\
\Df_{0} &:=& \Big\{(h_t)_{t\in K} \Big| h_t\in D_0~ \forall t,\ \sup_{s,t} d_D(h_t,h_s) < \infty\Big\}
\\
\Cf_{0} &:=& \Big\{(h_t)_{t\in K} \Big| h_t\in C_0~ \forall t,\ \sup_{s,t} d_D(h_t,h_s) < \infty\Big\}
\end{eqnarray*}
On the sets $\Rf_\Phi$ and $\Df_{0}$, define the metrics
\[
d_{R,\Phi}((h_t)_{t\in K},(g_t)_{t\in K}) := \sup_t d_R(h_t,g_t) \text{ and } d_{D,\Phi}((h_t)_{t\in K},(g_t)_{t\in K}) := \sup_t d_D(h_t,g_t),
\]
respectively. For elements $(h_t)_{t\in K}, (g_t)_{t\in K}$ set $(h_t)_{t\in K} + a (g_t)_{t\in K} := (h_t+ag_t)_{t\in K}$ and assume that with this definition, $(\Df_0,d_{D,\Phi})$ and $(\Rf,d_{R,\Phi})$ are metrized topological vector spaces.
Define the map
\[
\Phi:
\left\{
\begin{array}{ccc}
\Df_\Phi &\to& \Rf_\Phi
\\
(h_t)_{t\in K} &\mapsto& (\phi(h_t))_{t\in K}.
\end{array}
\right.
\]
\end{itemize}

\begin{theorem} \label{th:qhd}
Under setup (S) the map $\Phi$ is pseudo-Hadamard differentiable at $X:=(x)_{t\in K}$ tangentially to $\Uf\langle \Df_0\rangle$ where
\[
\Uf :=\Big\{ (h_t)_{t\in K}: h_t \in \Cf_0 ~ \forall t \in K, \sup_{d_K(s,t)\leq \delta} d_D(h_s,h_t) = o(1) \text{ as } \delta\to 0\Big\}
\]
and the derivative is given by
\[
\Phi_X':
\left\{
\begin{array}{ccc}
\Cf_0 &\to& \Rf
\\
(h_t)_{t\in K} &\mapsto& (\phi_x' h_t)_{t\in K}.
\end{array}
\right.
\]
\end{theorem}

Since quasi-Hadamard differentiability also implies classical Hadamard differentiability, the above result continues to hold in the classical setting.\\

As an illustration of the above result, let us consider the specific setting of Remark \ref{rem:gencom} in Section \ref{sec:genana} where $(R,d_R) = (\ell^\infty(\Gc_1)\times...\times\ell^\infty(\Gc_L),\|\cdot\|_{\max})$, $D = \ell^\infty(\Fc_1)\times...\times\ell^\infty(\Fc_J)$, $d_D((f_1,...,f_J),(f_1',...,f_J')) := \max_j \|f_j-f_j'\|_\infty$.
Consider the map
\[
\Psi_K:
\left\{
\begin{array}{ccc}
\Df_\Psi &\to& \Rf_\Psi
\\
(h_\kappa)_{\kappa\in K} &\mapsto& \Big(w(\kappa)\phi\Big(\frac{h_\kappa}{w(\kappa)}\Big)\Big)_{\kappa\in K}
\end{array}
\right.
\]
where $K$ is a compact set, $\Rf_\Psi \subset \Lc^\infty(\Gc_1,...,\Gc_L;K)$ and
\[
\Df_\Psi :=\Big\{(H_\kappa)_{\kappa\in K}\Big| \frac{H_\kappa}{w(\kappa)} \in D_\phi \forall \kappa\in K, \sup_{\kappa\in K}\Big\|w(\kappa)\phi\Big(\frac{H_\kappa}{w(\kappa)}\Big)\Big\|<\infty\Big\}.
\]
Note that with this definition, $\Vb_n = \alpha_n \Big(\Psi_K(w(\kappa)\hat y_{n,\kappa}) - \Psi_K(X_K)\Big)$ where $X_K := ((\kappa,f_1,...,f_J)\mapsto w(\kappa)x(f_1,...,f_J))$ and that $\Yb_n = \alpha_nw(\kappa)(\hat y_{n,\kappa} - X_K)$. As long as $\inf_{\kappa\in K} |w(\kappa)|>0$, compact differentiability of $\Psi_K$ with derivative
\[
(\Psi_K)_X':
\left\{
\begin{array}{ccc}
\Vf_K &\to& \Rf
\\
(h_\kappa)_{\kappa\in K} &\mapsto& (\phi_x' h_\kappa)_{\kappa\in K}.
\end{array}
\right.
\]
is a direct consequence of Theorem \ref{th:qhd} [here, $\Vf_K$ is defined similarly to $\Uf$ with $U$ replaced by $V$]. To see this, consider a sequence of real numbers $r_n \searrow 0$ and $h_n \in \Df$ such that $X_K + r_n h_n \in \Df_\Psi$ for all $n\in \N$ with $h_n\to h \in \Vf_K$. Then, by compact differentiability of $\Phi$,
\[
r_n^{-1}(\Psi_K(X_K + r_nh_n) - \Psi_K(X_K)) = r_n^{-1}w(\cdot)(\Phi_K(\tilde X_K + r_n\tilde h_n) - \Phi_K(\tilde X_K)) \to w(\cdot)\Phi_K'\tilde h.
\]
where $\tilde X_K := ((\kappa,f_1,...,f_J)\mapsto x(f_1,...,f_J))$, $\tilde h_n := ((\kappa,f_1,...,f_J)\mapsto h_n(\kappa,f_1,...,f_J)/w(\kappa))$ and $\tilde h:= ((\kappa,f_1,...,f_J)\mapsto h(\kappa,f_1,...,f_J)/w(\kappa))$. Finally, observing that $\Phi_K'$ is linear, compact differentiability of $\Psi_K$ and the definition of its derivative follow.\\

This result is of independent interest. For example, \cite{gazh2011} recently demonstrated that compact differentiability can be used to establish large and moderate deviation principles. The findings above allow to carry their results into the setting of statistics from subsamples and could for example be used to analyze rejection probabilities of various breakpoint tests.

\section{Applications}\label{sec:app}

In this section, we demonstrate how the results in Section \ref{sec:gen} can be applied to various subsample based methodologies studied in the recent literature. Throughout this section, we will assume that we have a sample of data $X_1,...,X_n$ from a strictly stationary time series. The process $\hat y_{n,s,t}$ is assumed to be based on the sub-sample $X_{\lfrf{ns}+1},...,X_{\lfrf{nt}}$, i.e. of the form given in equation \eqref{eq:xnt} .
In what follows, write $\theta = \phi(x)$ for the parameter of interest and define $\hat \vartheta_{n,s,t} := \phi(\hat y_{n,s,t})$. For notational convenience, we also consider the quantity $\hat\theta_{n,k,j}$ which is computed from the data $X_k,X_{k+1},...,X_j$. Note that $\hat\theta_{n,k,j} = \hat \vartheta_{n,k/n-1,j/n}$.
For the sake of a shorter notation, introduce the abbreviation $\Vb_{s,t} := \Vb(s,t,\cdot)$.

\subsection{Self-normalization}\label{sec:appsn}

For a weakly dependent stationary time series, inference on a finite-dimensional quantity (say, mean or median) typically involves
a consistent estimation of the asymptotic variance matrix of the sample estimator. The difficulty with the traditional approach lies in the
bandwidth parameter(s) involved in the consistent estimation, which also occurs for other existing approaches, such as sub-sampling [\cite{PR1994}], moving
block bootstrap [\cite{Kuns1989}] and block-wise empirical likelihood [\cite{Kita1997}].  To avoid the bandwidth selection, a general self-normalized approach to confidence interval construction and hypothesis testing for a stationary time series has been developed in \cite{shao2010}. The basic idea is to use recursive estimates to form an inconsistent estimator of asymptotic variance (matrix) of a statistic and use a non-standard but pivotal limiting distribution to perform the inference. The SN approach is convenient to implement as recursive estimates can be easily calculated with no need to develop new algorithms. Moreover, it does not involve any bandwidth parameters and its finite sample performance is comparable or could be  superior to some other existing bandwidth-dependent inference methods, as shown in \cite{shao2010}. Owing to these nice features, it has been recently extended to a few important inference problems in time series; see \cite{shzh2010,shao2011,shao2012,zhsh2013}, among others.

The theory for the SN approach was first developed in \cite{shao2010,shao2010c} by adopting a traditional approach, which is based on a linearization of the statistic and assumptions on uniform negligibility of remainder terms. More precisely, \cite{shao2010} assumed that
\begin{equation}\label{eq:lin}
\hat \theta_{n,1,k} = \theta + k^{-1}\sum_{i=1}^{k} L(X_i) + R_n(k/n)
\end{equation}
where $\{R_n(k/n)\}_{k=1}^{n}$ denote  negligible remainder terms. To describe the basic idea of Shao's approach, note that we generally expect that for a weakly dependent stationary time series and smooth functional $\phi$, $R_n(1) = o_P(n^{-1/2})$ and
\begin{eqnarray}
\label{eq:clt}
n^{-1/2}\sum_{j=1}^{n}L(X_j)\weak N(0,\Sigma),
\end{eqnarray}
where $\Sigma=\sum_{k\in\Z}\cov(L(X_0),L(X_k))>0$ is the so-called long run variance matrix.
Further note that we implicitly assume $\E[L(X_j)]=0$, which is trivially satisfied in many cases. Inference on $\theta$ is then based on estimating the covariance matrix $\Sigma$ consistently, which can be difficult as it involves a choice of bandwidth parameters. To avoid those complications, \cite{shao2010} proposed to consider the self-normalized statistic
\[
G_n = n(\hat\theta_{n,1,n}-\theta)'V_n^{-1}(\hat\theta_{n,1,n}-\theta)
\]
where $V_n=n^{-2}\sum_{j=1}^{n}j^2(\hat\theta_{n,1,j}-\hat{\theta}_{n,1,n})(\hat{\theta}_{n,1,j}-\hat{\theta}_{n,1,n})'$ is the self-normalization matrix. In Shao (2010a,b), the asymptotic distribution of $G_n$ was derived under the following assumptions:
\begin{eqnarray}
\label{eq:fclt}
n^{-1/2}\sum_{j=1}^{\lfloor nt\rfloor} L(X_j)\weak \Sigma^{1/2} \Bb(t),
\end{eqnarray}
 \begin{eqnarray}
\label{eq:neg2}
R_n(1) = o_P(n^{-1/2}),
\quad n^{-2}\sum_{j=1}^{n}|jR_n(j/n)|^2=o_p(1)
\end{eqnarray}
with $\Bb$ denoting a vector of independent Brownian motions on $[0,1]$. To verify (\ref{eq:neg2}), a common approach is to derive a uniform Bahadur representation
for $\hat{\theta}_{n,1,\lfloor nt\rfloor}$ and control the order of $R_n(t)$ uniformly over $t\in [0,1]$. Such a task is in general not easy and it requires a tedious case-by-case study. Under the assumptions above, \cite{shao2010} proved that
\begin{equation}\label{eq:sn}
G_n \weak \Bb^T(1)\Big(\int_0^1\Big(\Bb(t) - t\Bb(1)\Big)\Big(\Bb(t) - t\Bb(1)\Big)^T dt \Big)^{-1}\Bb(1),
\end{equation}
where the limiting distribution is pivotal and does not depend on the unknown covariance matrix $\Sigma$.

Using the results in Section \ref{sec:gen}, we can both considerably generalize the findings in \cite{shao2010} and at the same time avoid tedious calculations required to bound remainder terms. The key observation is that the only result required to derive \eqref{eq:sn} is weak convergence of the process $\Big(\sqrt{n}t(\hat \vartheta_{n,0,t} - \theta)\Big)_{t\in [0,1]}.$ In the language of Section \ref{sec:genana}, this amounts to setting $K = \{0\}\times[0,1]$. Assuming that $\phi(x)$ is an element of $\R^p$, the quantity $\Vb_n(s,t,\cdot)$ can be viewed as a $\R^p-$valued vector. Abusing notation, denote this vector by $\Vb_n(s,t)$. Similarly, denote by $\Vb_{s,t}$ the vector $\Vb(s,t,\cdot)$. Some straightforward calculations show that under assumptions (A1)-(A3), (C), (W) the statistic $G_n$ can be represented as
\[
G_n = \Vb_n^T(0,1)\Big(\int_0^1\Big(\Vb_n(0,t) - t\Vb_n(0,1)\Big)\Big(\Vb_n(0,t) - t\Vb_n(0,1)\Big)^T dt \Big)^{-1}\Vb_n(0,1) + o_P(1).
\]
An application of Theorem \ref{theorem:main1} with the set $K = \{0\}\times [0,1]$ in combination with the discussion at the beginning of this section and the continuous mapping theorem yields
\[
G_n \weak \Vb^T_{0,1}\Big(\int_0^1\Big(\Vb_{0,t} - t\Vb_{0,1}\Big)\Big(\Vb_{0,t} - t\Vb_{0,1}\Big)^T dt \Big)^{-1}\Vb_{0,1}.
\]
Under the assumption that $\Vb(0,t,\cdot) = \Sigma^{1/2}\Bb(0,t,\cdot)$, the limit of the statistic $G_n$ is pivotal. Note that the limiting process will typically have this form in most settings with weakly dependent data, see Remark \ref{rem:a1weak}.

With the general machinery of Section \ref{sec:gen} at hand, there are several extensions and remarks that can be made to the self-normalization approach. First, observe that we can replace the self-normalization matrix $V_n$ with a more general statistic of the form
\begin{eqnarray*}
V_n(H) &:=& \int_\Delta (\hat \vartheta_{n,s,t}-(t-s)\hat \vartheta_{n,0,1})(\hat \vartheta_{n,s,t}-(t-s)\hat \vartheta_{n,0,1})^T dH(s,t)
\end{eqnarray*}
with $H$ denoting an arbitrary probability measure on $\Delta$. By the continuous mapping theorem, we have joint convergence of $(V_n(H),\hat\vartheta_{n,0,1})$ to $(W(H),\Vb_{0,1})$ where
\[
W(H) := \int_\Delta \Big(\Vb_{s,t} - (t-s)\Vb_{0,1}\Big)\Big(\Vb_{s,t} - (t-s)\Vb_{0,1}\Big)^T dH(s,t).
\]
Assuming that $W(H)$ is non-singular almost surely [which happens as soon as $H$ places mass on sufficiently many different points], the asymptotic distribution of the generalized self-normalized statistic $G_n(H)$ follows. We thus have derived the following result.

\begin{prop}
Let assumptions (A1), (A2), (W), (C) hold and assume that either the support of $H$ is bounded away from the set $\{(t,t)|t\in[0,1]\}$ or that (A3) holds.
Additionally, assume that $W(H)$ is non-singular almost surely. Then the generalized SN-statistic $G_n$ satisfies
\[
G_n(H) := \Vb_n^T(0,1) V_n(H)^{-1}\Vb_n(0,1) \weak \Vb_{0,1}^T W(H)^{-1}\Vb_{0,1}.
\]
\end{prop}

Finally, note that by the discussion in Remark \ref{rem:clip} it might be advantageous to exclude estimators $\hat\theta_{n,k,l}$ that are based on too small proportions of data. By considering a modified version of the statistic $G_n$ of the form
$
\bar G_n(H) := \Vb_n^T(0,1) \bar V_n(H)^{-1}\Vb_n(0,1)
$
with
\[
\bar V_n(H) := \frac{\int_\Delta (\hat \vartheta_{n,s,t}-(t-s)\hat \vartheta_{n,0,1})(\hat \vartheta_{n,s,t}-(t-s)\hat \vartheta_{n,0,1})^T I\{t-s>n^{-\gamma}\}dH(s,t)}{\int_\Delta I\{t-s>n^{-\gamma}\}dH(s,t)},
\]
and arbitrary $\gamma\in(0,1/2)$, we would obtain the convergence $\bar G_n(H) \weak \Vb_{0,1}^T W(H)^{-1}\Vb_{0,1}^T$ without the need for assumption (A3).\\


\subsection{Subsampling and fixed-b corrections}\label{sec:appfb}

Sub-sampling [\cite{PR1994}] has been used in a wide range of inference problems for time series. The basic idea is that the distribution of an estimator computed from a sufficiently large sub-sample of the data should be close to that of the estimator from the whole data set. Confidence intervals and tests can then be constructed by approximating the unknown distribution of the estimator with sub-sampling counterparts. To accommodate the time series dependence non-parametrically, it involves the  sub-sampling window width $l$, which needs to go to infinity as sample size goes to infinity but at a slower rate to achieve consistent approximation. In practice, the choice of $l$ affects the sub-sampling distribution estimator and related operating characteristics, although its role does not show up in the conventional first order asymptotics. In \cite{shpo2013}, the traditional sub-sampling method was calibrated using a p-value based argument under the so-called fixed-$b$ asymptotics [\cite{KV2005}], where $b=l/n$. For simplicity, assume that $\theta$ is $\R^d$-valued. Defining $N = n - l + 1$, the sub-sampling based estimator of the distribution function of $\|\sqrt{n}\{\hat{\theta}_{n,1,n}-\theta\}\|$ evaluated at $x$ is
\[
L_{n,l}(x)=N^{-1}\sum_{j=1}^{N}I\{\|\sqrt{l}(\hat{\theta}_{n,j,j+l-1}-\hat{\theta}_{n,1,n})\|\le x\}.
\]
The corresponding p-value of the test statistic $\|\sqrt{n}(\hat{\theta}_{n,1,n} - \theta_0)\|$ for the null hypothesis $\theta=\theta_0$ is
\[
\hat p_n(b) =N^{-1}\sum_{j=1}^{N} I\{\|\sqrt{n}(\hat{\theta}_{n,1,n}-\theta_0)\|\le \|\sqrt{l}(\hat{\theta}_{n,j,j+l-1}-\hat{\theta}_{n,1,n})\|\}.
\]
Note that under the conditions $l/n+1/l = o(1)$ and additional regularity assumptions, $\hat p_n(b)$ has a uniform asymptotic distribution, see \cite{PRW1999}. Under the fixed-b asymptotic framework, $l/n=b\in (0,1]$ is held fixed. Following an elementary approach, the limiting null distribution of $\hat p_n(b)$, which equals
\[
G(b)=(1-b)^{-1}\int_0^{1-b} I\{\|\Sigma^{1/2}\Bb(1)\|\le \|\Sigma^{1/2}(\Bb(b+t)-\Bb(t)-b\Bb(1))\|/\sqrt{b}\}dt
\]
was derived in Shao and Politis (2013) by assuming that
\[
\hat \theta_{n,j,j+l-1} = \theta + l^{-1}\sum_{i=j}^{j+l-1} L(X_i) + R_n(j,j+l-1),
\]
that a similar representation holds for $\hat \theta_{n,1,n}$, that (\ref{eq:fclt}) holds for $\{L(X_t)\}$ with remainder $R_n(1,n)$, and that the remainder terms satisfy $\sqrt{n}|R_n(1,n)|=o_p(1)$ and $\sqrt{l}\sup_{j=1,\cdots,N}|R_n(j,j+l-1)|=o_p(1)$. Verifying the latter assumption for general functionals can be quite tedious and challenging. \\
Now, consider the general setup of Section \ref{sec:gen} and let conditions (C), (W), (A1) and (A2) hold. We apply Theorem \ref{theorem:main1} with $K := \{(t,t+b)|t\in [0,1-b]\}\cup\{(0,1)\}$
 and assume that the map
 $$h \mapsto \frac{1}{1-b}\int_0^{1-b} I\{ \|h(0,1)\| \leq \|h(t,t+b)-bh(0,1)\|/\sqrt{b} \}dt$$
 is continuous on a set of functions that contains the sample paths of $\|\Vb\|$ with probability one. In particular, this is the case if $\Vb(s,t) = (t-s)\Sigma^{1/2}(\Bb(t) - \Bb(s))$ with $\Sigma$ denoting a non-singular matrix and $\Bb$ a vector of independent Brownian motions [see the arguments in \cite{shpo2013}], which is  typically the case for  weakly dependent stationary time series. From now on, assume that this is the case. Observe that for $\theta = \theta_0$ we have in the setting discussed above
\[
\hat p_n(b) = \frac{1}{1-b}\int_0^{1-b} I\{\Vb_n(0,1) \leq \|\Vb_n(t,t+b)-b\Vb_n(0,1)\|/\sqrt{b} \}dt + o_P(1),
\]
where the negligibility of remainder follows from an application of the continuous mapping theorem. The results in Theorem \ref{theorem:main1} in combination with the continuous mapping theorem thus yield
\[
\hat p_n(b) \weak \frac{1}{1-b}\int_0^{1-b} I\{\|\Vb(0,1,\cdot)\| \leq \|\Vb(t,t+b,\cdot) - b\Vb(0,1,\cdot)\|/\sqrt{b} \} dt := P
\]
as soon as assumptions (C), (W), (A1) and (A2) hold.\\
Unless $\theta$ is real-valued, the asymptotic distribution of the statistic $\hat p_n(b)$ is in general not pivotal. \cite{shpo2013} proposed to estimate its distribution based on further sub-sampling. An alternative is to consider block bootstrap approximations such as those discussed in Section \ref{sec:genpro}. More precisely, consider a bootstrap version for $y_{n,s,t}$ which is of the form given in \eqref{eq:ynb1} and denote it by $\hat y_{n,s,t}^B$. Define a bootstrap version for $\hat\vartheta_{n,s,t}$ through $\hat\vartheta_{n,s,t}^B := \phi(\hat y_{n,s,t}^B)$. Assume that the map $\phi$ is continuous. Now Theorem \ref{th:boot1} combined with the continuous mapping theorem for the bootstrap in probability [see Theorem 10.8 in \cite{kosorok2008}] directly yields that under condition (WB)
\[
\hat p_n^B(b) := \frac{1}{1-b}\int_0^{1-b} I\{\sqrt{n}\|\hat\vartheta_{n,0,1}^B - \theta\| \leq \sqrt{nb}\|\hat\vartheta_{n,t,t+b}^B-\hat\vartheta_{n,0,1}^B\| \}dt \weakPi{M} P.
\]
Finally, note that the reasoning above does not rely on $\theta$ being $\R^p$-valued and that it is thus also possible to handle infinite dimensional parameters.

\subsection{Testing for change points}\label{sec:appcp}

Testing change points in a time series is a well-studied topic in econometrics and statistics; see \cite{P2006} for a recent review.
A large class of tests in the literature is based on the so-called CUSUM (cumulative sum) process and
the test statistic is a smooth functional of the CUSUM  process with Kolmogorov-Smirnov ($L^{\infty}$) test
and Cramer-von-Mises ($L^2$) test being two prominent examples.
 To accommodate the time series dependence and make the limiting null distribution pivotal, one needs to obtain a consistent estimator of the long run variance
 as a studentizer. As mentioned previously, consistent estimation involves a bandwidth parameter, the choice of which is even more difficult in the change point testing problem. In particular, the fixed bandwidth (e.g., $n^{1/3}$) is not adaptive to the magnitude of dependence and the data-dependent bandwidth could lead to the so-called non-monotonic power problem [\cite{Vo1999}], i.e., the power of the test can decrease when the alternative gets farther away from the null. To overcome the non-monotonic power problem, \cite{shzh2010} proposed SN-based tests in a general framework. Let $\theta_t = T(\mathcal{D}(X_t))\in \R^q$ be the quantity of interest which depends on the distribution of $X_t$ denoted by $\mathcal{D}(X_t)$. The goal is to test if there is a change point
in $\{\theta_t\}_{t=1}^{n}$, i.e.
\[H_0: \theta_1=\cdots=\theta_n\]
and the alternative hypothesis is
\[H_1: \theta_1=\cdots=\theta_{k^*}\not=\theta_{k^*+1}=\cdots=\theta_n~\mbox{for some unknown}~k^*,~1\le k^*< N.\]
This framework is general enough to include mean, median, autocorrelation at certain lags of a univariate time series. Let $T_n(k)=k/\sqrt{n}(\hat{\theta}_{n,1,k}-\hat{\theta}_{n,1,n})$ and
\begin{eqnarray*}
V_n(k) &=& n^{-2}\left\{\sum_{t=1}^{k} t^2 (\hat{\theta}_{n,1,t}-\hat{\theta}_{n,1,k}) (\hat{\theta}_{n,1,t}-\hat{\theta}_{n,1,k})'\right.
\\
&& \left.+\sum_{t=k+1}^{n}(n-t+1)^2 (\hat{\theta}_{n,t,n}-\hat{\theta}_{n,k+1,n})(\hat{\theta}_{n,t,n}-\hat{\theta}_{n,k+1,n})'\right\}.
\end{eqnarray*}
Then the test statistic is defined as $G_n=\sup_{k=1,\cdots,n-1} T_n(k)'V_n(k)^{-1}T_n(k)$. The asymptotic null distribution was derived in \cite{shzh2010} using an elementary approach. Specifically, they rely on the expansion of $\hat{\theta}_{n,t_1,t_2}$, i.e.,
\[
\hat{\theta}_{n,t_1,t_2}=\theta + (t_2-t_1+1)^{-1}\sum_{t=t_1}^{t_2} L(X_t) + R_n(t_1,t_2).
\]
Again the functional central limit theorem is assumed for $\{L(X_t)\}$ (i.e., (\ref{eq:fclt}) holds) and the remainder terms are assumed to be asymptotically negligible. In particular, \cite{shzh2010} assume that
\begin{eqnarray}
\label{eq:neg3}
\sup_{k=1,\cdots,n} |k R_n(1,k)|=o_p(n^{1/2})\quad \mbox{and} \quad \sup_{k=1,\cdots,n}|k R_n(n-k+1,n)|=o_p(n^{1/2}).
\end{eqnarray}
The above condition~(\ref{eq:neg3}) is not easy to verify and a detailed case-by-case study is needed.\\
Alternatively, consider the setting of Section \ref{sec:gen}. Under conditions (W), (A1), (A2) and (A3) with $K = (\{0\}\times[0,1])\cup([0,1]\times\{0\})$ it is possible to show that $G_n = o_P(1) + \sup_{r\in [0,1]} H_n(r)$ where
\[
H_n(r) := (\hat\vartheta_{n,0,r}-\hat\vartheta_{n,0,1})^T \hat W_{r,n}^{-1} (\hat\vartheta_{n,0,r}-\hat\vartheta_{n,0,1})
\]
with
\[
\hat W_{r,n} :=  \int_0^r (\hat\vartheta_{n,0,s}-\hat\vartheta_{n,0,r})^T(\hat\vartheta_{n,0,s}-\hat\vartheta_{n,0,r}) ds + \int_r^1 (\hat\vartheta_{n,s,1}-\hat\vartheta_{n,r,1})^T(\hat\vartheta_{n,s,1}-\hat\vartheta_{n,r,1})ds.
\]
Applying Theorem \ref{theorem:main1} in combination with the continuous mapping theorem yields weak convergence of $G_n$ to
\[
\sup_{r\in [0,1]} (\Vb_{0,r}-r\Vb_{0,1})^T W_r^{-1} (\Vb_{0,r}-r\Vb_{0,1})
\]
where
\[
W_r := \int_0^r \Big(\Vb_{0,s}-\frac{s}{r}\Vb_{0,r}\Big)^T\Big(\Vb_{0,s}-\frac{s}{r}\Vb_{0,r}\Big) ds + \int_r^1 \Big(\Vb_{s,1}-\frac{1-s}{1-r}\Vb_{r,1}\Big)^T\Big(\Vb_{s,1}-\frac{1-s}{1-r}\Vb_{r,1}\Big)ds.
\]

Finally, note that by considering the modification $\tilde G_n := \sup_{r\in [n^{-\gamma},1-n^{-\gamma}]} H_n(r)$ with $\gamma\in(0,1/2)$ arbitrary, assumption (A3) can be dropped. See Remark \ref{rem:clip} for further details.

\begin{appendix}
\section{Proofs of main results}

\textbf{Proof of Theorem \ref{theorem:main1}}
The proof consists of two steps. First, we show that the convergence holds for $K$ with $\inf_{(s,t)\in K} t-s > 0$, and second, we extend the result to general sets $K\subset \Delta$ under assumption (A3). The first step follows by an application of the functional delta method [see Theorem 3.9.4 in \cite{vandervaart1996}] in combination with Theorem \ref{th:qhd}. In particular, the space $\Df_0$ can be identified with $\Lc(\Fc_1,...,\Fc_J;K)$ since a finite norm $\|H\|_\Lc$ is equivalent to the distance $\|H-(x)_{(s,t)\in K}\|_\Lc$ being finite. Similarly, the space $R_\Psi$ is identified with $\ell^\infty(\Gc \times K)$ and the metric $d_{R,\Phi}$ corresponds to the supremum norm on $\Lc^\infty(\Gc_1,...,\Gc_L;K)$. Observe that
\[
\alpha_n(\hat y_{n,s,t}(f_1,...,f_J) - x(f_1,...,f_J)) = \frac{1}{t-s}~ \Yb_n(s,t,f_1,...,f_J).
\]
The functional delta method in combination with elementary considerations thus implies
\[
\alpha_n(\phi(\hat y_{n,s,t}) - \phi(x)) \weak \frac{1}{t-s}\phi_x'\Yb(s,t,\cdot) \quad \text{in} \quad \Lc^\infty(\Gc_1,...,\Gc_L;K),
\]
the factor $\frac{1}{t-s}$ can be moved in front since $\phi_x'$ is a linear map. Multiplying both sides by $t-s$, the Continuous Mapping Theorem [see Theorem 1.3.6 in \cite{vandervaart1996}] completes the first step of the proof.\\
For the second step, define the set $K_{S} := \{(s,t)\in K| t-s\in S\}$ and consider the approximating processes
\[
A_{i,n} := (t-s)\alpha_n I_{(s,t)\in K_{[1/i,1]}}(\phi(\hat y_{n,s,t}) - \phi(x)), \quad A_{i} := \phi_x'\Yb(s,t,\cdot)I_{(s,t)\in K_{[1/i,1]}}.
\]
It then suffices to verify the following three statements [see \cite{budevo2011}]
\begin{eqnarray*}
& (\rm{i}) & ~ \text{For every } i\in\N: ~ A_{i,n} \weak A_{i} ~\text{ for } {n\rightarrow\infty}, \\
& (\rm{ii}) & ~ A_{i} \weak \phi_x'\Yb(s,t,\cdot)~ \text{ for } {i\rightarrow\infty}, \\
& (\rm{iii}) & ~ \text{For every }\eps>0: ~ \lim_{i\rightarrow\infty}\limsup_{n\rightarrow\infty} \Pb^*(\|A_{i,n}-\Vb_n\|>\varepsilon) = 0.
\end{eqnarray*}
The first statement is the weak convergence established in the first step. For (ii), note that
\[
\Big\|A_{i} - \phi_x'\Yb(s,t,\cdot)\Big\| \leq \|\phi_x'\|_{op}\sup_{(s,t) \in K_{[0,1/i]}}\sup_j\sup_{f\in\Fc_j}|\Yb_j(s,t,f)|,
\]
[here, $\|\cdot\|_{op}$ denotes the operator norm] and the right-hand side converges to zero in probability, this is a direct consequence of assumption (A1).

Finally, for a proof of (iii) note that for $\beta_n := \gamma_n\vee \alpha_n^{-1/2} = o(1)$ from Lemma \ref{lem:tech1}
\begin{eqnarray*}
\|A_{i,n} - \Vb_n\| &\leq& \alpha_n\sup_{(s,t)\in K_{[\beta_n n^{-1/2},1/i]}}(t-s)\|\phi(\hat y_{n,s,t}) - \phi(x)\| + \sup_{(s,t)\in K_{[0,\beta_n n^{-1/2}]}}(t-s)\|\phi(\hat y_{n,s,t}) - \phi(x)\|
\\
&\leq&  C\sup_{(s,t)\in K_{[\beta_n n^{-1/2},1/i]}}\sup_j\sup_{f\in\Fc_j}|\Yb_{n,j}(s,t,f)| + \sup_{(s,t)\in K_{[0,\beta_n n^{-1/2}]}}(t-s)\Big(\|\phi(\hat x_{n,t})\| + \|\phi(x)\|\Big)
\\
&& + I\Big\{\sup_{(s,t) \in K_{[\beta_n n^{-1/2},1/i]}}(t-s)\alpha_n\|\hat y_{n,s,t}-x\| > \eps\Big\}\alpha_n\sup_{(s,t)\in K_{[\beta_n n^{-1/2},1/i]}}\|\phi(\hat y_{n,s,t}) - \phi(x)\|
\\
&=:& R_{n,1} + R_{n,2} + R_{n,3}.
\end{eqnarray*}
Here the second inequality follows by an application of Lemma \ref{lem:phi} on the set $\Big\{\sup_{(s,t)\in K_{[\beta_n n^{-1/2},1/i]}}\alpha_n\|\hat y_{n,s,t}-x\| \leq \eps\Big\}$ after observing that by definition
\[
\alpha_n\|\hat y_{n,s,t}-x\| = \frac{1}{t-s}\sup_j\sup_{f\in\Fc_j}|\Yb_{n,j}(s,t,f)|.
\]
Condition (A3) implies that $R_{n,2} = o_P^*(1)$.
To see that $R_{n,1}+R_{n,3}$ converge to zero in outer probability, define the set
\[
S_j(i,\eps) := \Big\{y \in \ell^\infty(K\times \Fc_j)\Big| \sup_{(s,t) \in K_{[0,1/i]}}\sup_{f\in\Fc_j}|y((s,t),f)| \geq \eps \Big\}.
\]
This set is closed, and by the Portmanteau theorem [Theorem 1.3.4 in \cite{vandervaart1996}] combined with the weak convergence of $\Yb_{nj}$ and assumption (A1) on $\Yb$ we obtain
\[
\limsup_{n\to\infty}P^*(\Yb_{nj} \in S_j(i,\eps)) \leq P(\Yb_{j} \in S_j(i,\eps))
\]
for $j=1,...,J$. By condition (A1), $\lim_{i\to \infty}  P(\Yb_{j} \in S_j(i,\eps)) = 0$ for every $\eps>0$.

This shows that $R_{n,1} = o_P^*(1)$ and $R_{n,3} = o_P^*(1)$. Thus the proof is complete. \hfill $\Box$\\
\\
\textbf{Proof of Theorem \ref{th:boot1}} The first assertion follows by an application of the bootstrap functional delta method [see e.g. Theorem 12.1 in \cite{kosorok2008}]. For more details on the appropriate identification of spaces, see the proof of Theorem \ref{theorem:main1} in the present note. In order to prove the second part, define the set $K_{S} := \{(s,t)\in K| t-s\in S\}$ and consider
\[
A_{i,n}^b := (t-s)\alpha_n I_{(s,t) \in K_{[1/i,1]}}(\phi(\hat y_{n,s,t}^b) - \phi(\hat y_{n,s,t})), \quad A_{i} := \phi_x'\Yb(t,\cdot)I_{(s,t) \in K_{[1/i,1]}}.
\]
By Lemma \ref{lem:bootappr} it then suffices to verify the following three statements which can be regarded as adaptation of Theorem 4.2 in \cite{billingsley1968} to the present setting
\begin{eqnarray*}
& (\rm{i}) & ~ \text{For every } i\in\N: ~ A_{i,n}^b \weakPi{M} A_{i} ~\text{ for } {n\rightarrow\infty}, \\
& (\rm{ii}) & ~ A_{i} \weak \phi_x'\Yb(s,t,\cdot)~ \text{ for } {i\rightarrow\infty}, \\
& (\rm{iii}) & ~ \text{For every }\varepsilon>0: ~ \lim_{i\rightarrow\infty}\limsup_{n\rightarrow\infty} \Pb^*(\|A_{i,n}^b-\Vb_n^b\|>\varepsilon) = 0.
\end{eqnarray*}
Assertion (i) follows from the first part. Assertion (ii) can be established by exactly the same arguments as the corresponding statement in the proof of Theorem \ref{theorem:main1}. For a proof of the third assertion, note that $\Yb_n^b \weakPi{M} \Yb$ implies $\Yb_n^b \weak \Yb$, see e.g. the proof of Theorem 10.4, assertion $(ii)\Rightarrow (i)$ in \cite{kosorok2008}. Thus assertion (iii) follows by exactly the same arguments as (iii) in the proof of Theorem \ref{theorem:main1}. Hence the proof is complete. \hfill $\Box$\\
\\
\textbf{Proof of Proposition \ref{prop:mult}}
Observe the representation
\[
\hat y_{n,s,t}(\cdot) - x(\cdot) = \frac{n}{(\lfrf{nt}-\lfrf{ns})\vee 1}\Big(\frac{\lfrf{nt}}{n}~(\hat x_{n,t} - x) - \frac{\lfrf{ns}}{n}~(\hat x_{n,s} - x)\Big)
\]
and thus setting $'0/0=0'$
\[
\Yb_{n}(s,t,\cdot) = \frac{n(t-s)}{(\lfrf{nt}-\lfrf{ns})\vee 1}\Big(\frac{\lfrf{nt}}{nt}~\Gb_n(t,\cdot) - \frac{\lfrf{ns}}{ns}~\Gb_n(s,\cdot)\Big)
\]
Observe that $\Big|\frac{\lfrf{nt}}{nt} - 1\Big|\leq \frac{1}{\lfrf{nt}\vee 1}$ and $\Big|\frac{n(t-s)}{(\lfrf{nt}-\lfrf{ns})\vee 1} -1 \Big|\leq \frac{3}{(\lfrf{nt}-\lfrf{ns})\vee 1}$. Defining $\tilde \Yb_n(s,t,\cdot) := \Gb_n(t,\cdot) - \Gb_n(s,\cdot)$, observe that $\tilde \Yb_n \weak \Yb$ by the continuous mapping theorem. Moreover, $\sup_{t} \Big\|\frac{\lfrf{nt}}{nt}~\Gb_n(t,\cdot) - \Gb_n(t,\cdot)\Big\| = o_P^*(1)$ since for $t\geq n^{-1/4}$ the factor $\frac{\lfrf{nt}}{nt}$ tends to one uniformly and since $\sup_{t\leq n^{-1/4}}\|\Gb_n(t,\cdot)\| = o_P^*(1)$ by arguments similar to those used to establish the negligibility of $R_{n,1}$ at the end of the proof of Theorem \ref{theorem:main1}. Thus it remains to show that $\Big(\frac{n(t-s)}{(\lfrf{nt}-\lfrf{ns})\vee 1} - 1\Big)\tilde \Yb_n(s,t,\cdot)$ is uniformly small. This can be done by similar arguments [distinguish the cases $t-s \leq n^{-1/4}$ and $t-s > n^{-1/4}$]. This completes the proof. \hfill $\Box$\\
\\
\textbf{Proof of Theorem~\ref{th:donsker}} Since it suffices to show asymptotic tightness of each process $\Gb_{n,j}$ individually, we will focus on $\Gb_{n,1}$. To simplify notation, define $\Z_n := \Gb_{n,1}$, $\Fc := \Fc_1$, $\Fc_\delta := \Fc_{1,\delta}$.
Start by noting that under the assumptions of the theorem together with \eqref{eq:res3} we have for some finite constant $C_1$
\begin{equation}\label{eq:res2}
\sup_{k\in\N}\E^* \|\Z_k(1,\cdot)\|^q_{\Fc} \le C_1 < \infty
\end{equation}
To see this, fix $\delta>0$ and cover the set $\Fc$ with $N$ balls of radius $\delta$ and centers $f_1,...,f_N$. Then make use of the bound
\[
\sup_{k\in\N}\E^* \|\Z_k(1,\cdot)\|^q_{\Fc} \leq \max_{1\leq k\leq n_0} \E^*\|\Z_k(1,\cdot)\|^q_{\Fc} + \max_{j=1,...,N}\sup_{k\in\N}\E^* \|\Z_k(1,f_j)\|^q + \sup_{n\geq n_0} \E^*\|\Z_n(1,\cdot)\|_{{\cal F}_{\delta}}^q
\]
and condition \eqref{eq:res3}.\\
In order to establish asymptotic tightness of $\Gb$, apply Theorem 1.5.7 in \cite{vandervaart1996} with the metric $d((s,f),(t,g)) := \rho(f,g)+|s-t|$. By the triangle inequality, we have
\[
\sup_{|s-t|+\rho(f,g)<\delta} |\Z_n(s,f)-\Z_n(t,g)|\le \sup_{0\le t\le 1}\|\Z_n(t,f)\|_{{\cal F}_{\delta}} +  \sup_{|s-t|<\delta} \|\Z_n(s,f)-\Z_n(t,f)\|_{\cal F}
\]
Start by considering the first term. Define $S_k(g)=\sum_{j=1}^{k}\{g(X_i)-\E g(X_i)\}$ and note that
\[
\sup_{0\le t\le 1}\|\Z_n(t,f)\|_{{\cal F}_{\delta}}
\leq \max_{1\leq k\leq n} \sqrt{\frac{k}{n}}\|\Z_k(1,\cdot)\|_{\Fc_
\delta}
= \frac{1}{\sqrt{n}}\max_{1\leq k\leq n}\|S_k\|_{\Fc_
\delta}.
\]
Fix $\epsilon\in (0,\{1-2^{-1/2+1/q}\}^{q/(q-1)}/2^{q/(2q-2)})$. Under (\ref{eq:res1}),  there exists a $\delta_0>0$ and $n_0\in\N$, such that when $\delta\in (0,\delta_0)$ and $n\ge n_0(\epsilon)$, $\E^*\|\Z_n(1,\cdot)\|_{{\cal F}_{\delta}}^q\le \epsilon^2$. Moreover, under \eqref{eq:res2} we have $\max_{1\leq k\leq n_0}\|S_k\|_{\Fc_\delta} \leq C_1\sqrt{n_0}$ for all $n_0 \in \N$.
By the Markov inequality and Proposition 1(ii) in \cite{Wu2007}, for $q>2$, $d=d_n=\lfloor \log n/(\log 2)\rfloor +1$,
\begin{eqnarray*}
P^*(\max_{1\le k\le n}\|S_k\|_{{\cal F}_{\delta}}>\sqrt{n}\epsilon)&\le& (\sqrt{n} \epsilon)^{-q}\E^*[\max_{1\le k\le n} \|S_k\|_{{\cal F}_{\delta}}^q ]\\
&\le&(\sqrt{n} \epsilon)^{-q}\left(\sum_{j=0}^{d}2^{(d-j)/q}\left\{\E^*\sup_{g\in {\cal F}_{\delta}} |S_{2^j}|^q\right\}^{1/q}\right)^{q}
\\
&\le& \epsilon^{-q}n^{-q/2}\left(O(n)+\left\{\sum_{j=\lfloor \log n_0/(\log 2)\rfloor+1}^{d}2^{(d-j)/q} (\epsilon^{2q} 2^{jq/2})^{1/q}\right\}^{q}\right)\\
&\le &\epsilon^{-q}O(n^{-q/2+1})+n^{-q/2}\epsilon^q\left(\sum_{j=\lfloor \log n_0/(\log 2)\rfloor+1}^{d}2^{(d-j)/q} (2^{jq/2})^{1/q}\right)^{q}\\
&\le&\epsilon^{-q}O(n^{-q/2+1})+\frac{2^{q/2}}{\{1-2^{-1/2+1/q}\}^q}\epsilon^q\\
&<&\epsilon
\end{eqnarray*}
for $n\ge n_1(\epsilon)\vee n_0(\epsilon)$, where $n_1(\epsilon)$ is chosen such that the last inequality holds. Since $\epsilon$ was arbitrary, we have shown that $\limsup_{n\to\infty}P^*(\sup_{0\le t\le 1}\|\Z_n(t,f)\|_{{\cal F}_{\delta}} >\epsilon)<\epsilon$ for all $\delta < \delta_0$. It thus remains to consider the second term.
Since the increments of $\Z_n(s,f)$ in $s$ are stationary, the above probability can be bounded by
\begin{eqnarray}
\label{eq:bound2}
\lceil \frac{1}{\delta} \rceil P^*\left(\sup_{0\le s\le \delta} \|\Z_n(s,f)\|_{\cal F}>\epsilon \right)=\lceil \frac{1}{\delta} \rceil P^*\left(\max_{1\le k\le n\delta} \|\sqrt{k}\Gb_k\|_{\cal F}>\sqrt{n}\epsilon\right).
\end{eqnarray}
Let $d(\delta)=\lfloor \log (n\delta)/(\log 2)\rfloor +1$. Again by the Markov inequality and Proposition 1(ii) in \cite{Wu2007},
\begin{eqnarray*}
P^*\left(\max_{1\le k\le n\delta} \|S_k\|_{\cal F}>\sqrt{n}\epsilon\right)&\le&(\sqrt{n}\epsilon)^{-q}\E^* \max_{1\le k\le n\delta} \|S_k\|_{\cal F}^q\\
&\le&(\sqrt{n}\epsilon)^{-q}\left\{\sum_{j=0}^{d(\delta)} 2^{(d(\delta)-j)/q} \left(\E^*\|S_{2^j}\|_{\cal F}^q\right)^{1/q}\right\}^{q}\\
&\le&(\sqrt{n}\epsilon)^{-q}\left\{\sum_{j=0}^{d(\delta)} 2^{(d(\delta)-j)/q} C_1^{1/q} 2^{j/2} \right\}^{q}\\
&\le&(\sqrt{n}\epsilon)^{-q} C_1 2^{d(\delta)q/2}\frac{1}{\{1-2^{-(1/2-1/q)}\}^q}\\
&\le&C_2 \epsilon^{-q} \delta^{q/2}
\end{eqnarray*}
for $n$ sufficiently large. Combined with~(\ref{eq:bound2}), we get $ \limsup_{n\to\infty}P^*(\max_{0\le j\delta\le 1} \sup_{j\delta\le s\le (j+1)\delta} \|\Z_n(s,f)-\Z_n(j\delta,f)\|_{\cal F}>\epsilon)<\epsilon$ when $\delta<(\epsilon^{q+1}/C_2)^{1/(q/2-1)}$. The proof is thus complete.
\qed\\
\\
\textbf{Proof of Theorem \ref{th:qhd}}
Let $a_n = o(1)$ and $H^{(n)}$ denote a sequence in $\Df_0$ with $H^{(n)} \to H \in \Uf$ such that $X+a_nH^{(n)} \in \Df_\Phi \ \forall n\in\N$. We need to show that
\[
a_n^{-1}\Big(\Phi(X+a_nH^{(n)}) - \Phi(X)\Big) \to \Phi_X'H.
\]
Assume that this does not hold. Then there exists a sequence $t_n$ and a positive number $b$ such that
\begin{equation} \label{eq2}
d_R\Big(a_n^{-1}(\phi(x + a_nH^{(n)}_{t_n}) - \phi(x)),(\phi_x'\cdot H_{t_n})\Big) \geq b
\end{equation}
for all $n \geq N_0$. On the other hand, the sequence $H^{(n)}_{t_n}$ has a subsequence $H^{(n_k)}_{t_{n_k}}$ which converges to $H_{t_\infty}$ for some $t_\infty \in K$. To see that this is the case, start by noting that $t_n$ is a sequence in a compact metric space, i.e. it has a convergent subsequence $t_{n_k} \to t_\infty$ with $t_\infty \in K$. The definition of the set $\Uf$ then implies that $H_{t_{n_k}} \to H_{t_\infty}$. Together with the uniform convergence $\sup_t d_D(H^{(n)}_t,H_t) = o(1)$ this yields $H^{(n_k)}_{t_{n_k}} \to H_{t_\infty}$. Now quasi compact differentiability of $\phi$ tangentially to $C_0\langle D_0\rangle$ implies
\[
a_n^{-1}\Big(\phi(x + a_nH^{(n_k)}_{t_{n_k}}) - \phi(x)\Big) \to \phi_x' H_{t_\infty},
\]
and together with continuity of $\phi_x'$ this contradicts (\ref{eq2}). Thus the proof is complete. \hfill $\Box$

\section{Auxiliary technical results}

\begin{lemma} \label{lem:phi}
Denote by $(R,\|\cdot\|_R)$ a normed vector space. Consider a second vector space $D$ with subsets $D_\phi,D_0 \subset D, C_0 \subset D_0$ and assume that $(D_0,\|\cdot\|_D)$ is a normed vector space. Let $\phi: D_\phi \to R$ be quasi compactly differentiable at $x$ tangentially to $C_0\langle D_0\rangle$ and assume $0 \in C_0$. Then there exist constants $\eps>0, K<\infty$
such that
\begin{equation} \label{eq1}
\|\phi(x)-\phi(x+y)\|_R \leq K\|y\|_D \quad \forall y\in D_0:\  \|y\|_D \leq \eps, x+y\in D_\phi.
\end{equation}
\end{lemma}
\textbf{Proof} Assume that (\ref{eq1}) does not hold. Then for any pair $\eps > 0, K<\infty$ there exists a $y_{K,\eps} \in D_0$ such that $x+y_{K,\eps} \in D_\phi$, $\|y_{K,\eps}\|_D\leq \eps$ and $\|\phi(x)-\phi(x+y_{K,\eps})\|_R > K\|y_{K,\eps}\|_D$. Consider the sequence $z_n := y_{n^{2},n^{-2}}$ and define $\alpha_n := \|z_n\|_D \neq 0$. Then
\[
\Big\| \frac{\phi(x +  n\alpha_n(n\alpha_n)^{-1}(z_n)) - \phi(x)}{n\alpha_n} \Big\|_R > \frac{n^2 \alpha_n}{n\alpha_n} = n \longrightarrow \infty.
\]
Moreover $\|(n\alpha_n)^{-1}(z_n)\|_D = n^{-1} = o(1)$, i.e. $(n\alpha_n)^{-1}(z_n) \rightarrow 0$.
This yields a contradiction since quasi compact differentiability of $\phi$ implies that [note that $n\alpha_n \leq n^{-1} = o(1)$]
\[
\frac{\phi(x +  n\alpha_n(n\alpha_n)^{-1}(z_n)) - \phi(x)}{n\alpha_n} \longrightarrow \phi_x' 0.
\]
Thus the proof is complete. $\hfill \Box$

\begin{lemma}\label{lem:tech1}
Under assumptions (W) and (A1) there exists a sequence of real numbers $\gamma_n = o(1)$ such that $\sup_{(s,t)\in K_{\alpha_n^{-1}\gamma_n}}\| \hat y_{n,s,t} - x \| = o_P^*(1)$ where we defined $K_{a} := \{(s,t)\in K| t-s \geq a\}$.
\end{lemma}
\textbf{Proof}
Define $K_{a}^C$ as the complement of $K_a$ in $K$ and set
\[
B_n := \sup_{(s,t)\in K_{\alpha_n^{-1/2}}^C} \sup_{f_1,...,f_J}\|\Yb_n(s,t,f_1,...,f_J)\|.
\]
By asymptotic equicontinuity of $\Yb_n$ [see the discussion in the proof of Theorem \ref{theorem:main1} for more details and note that $\sup_{s=t}\sup_{f_1,...,f_J}|\Yb_n(s,t,f_1,...,f_J)|\equiv 0$ a.s.] we have $B_n = o_P^*(1)$. This implies
\[
\forall\eps>0 \ \exists n_0(\eps)\in \N:\quad  (*)\  \forall n\geq n_0(\eps) \ P^*(B_n>\eps)<\eps.
\]
Note that $a \mapsto n_0(a)$ is decreasing since for any $a<b$ we have $P^*(B_n>a)<a \Rightarrow P^*(B_n>b)<b$. Set $N_0(\eps) := 2\inf\{n_0(\eps)| (*) \mbox{ holds}\}$ and define
\[
\delta_n := 2\inf\{\eps>0|n>N_0(\eps)\}.
\]
By construction $N_0(\delta_n) < n$, and thus $P^*(B_n>\delta_n)<\delta_n$. Moreover, $\delta_n \rightarrow 0$ since by construction $\delta_n \leq \eps \ \forall n \geq N_0(\eps/3)$. Defining $\gamma_n = \delta_n^{1/2}$ yields $B_n = o_P^*(\gamma_n)$. Note that
\[
\sup_{(s,t)\in K_{\alpha_n^{-1}\gamma_n}} \| \hat y_{n,s,t} - x\| \leq \sup_{(s,t)\in K_{\alpha_n^{-1/2}}} \|  \hat y_{n,s,t} - x \| + \sup_{(s,t) \in K_{\alpha_n^{-1/2}}^C\cap K_{\gamma_n\alpha_n^{-1}}} \| \hat y_{n,s,t} - x \|.
\]
Now observe that
\[
\sup_{(s,t)\in K_{\alpha_n^{-1/2}}} \| \hat y_{n,s,t} - x \|
\leq 2\alpha_n^{-1/2} \sup_{(s,t)\in K_{\alpha_n^{-1/2}}}\|\Yb_n(s,t,\cdot)\| = o_P^*(1),
\]
by arguments similar to those used to establish the negligibility of $R_{n,1}$ at the end of the proof of Theorem \ref{theorem:main1}. Similarly
\[
\sup_{(s,t) \in K_{\alpha_n^{-1/2}}^C\cap K_{\gamma_n\alpha_n^{-1}}} \| \hat y_{n,s,t} - x \| \leq \gamma_n^{-1} \sup_{(s,t) \in K_{\alpha_n^{-1/2}}^C}|\Yb_n(s,t,\cdot)\| = \gamma_n^{-1}B_n  = o_P^*(1)
\]
This completes the proof.
\hfill $\Box$

\begin{lemma}\label{lem:bootappr}
Given a sequence of random variables $M_1,M_2,...$, and a sequence of random elements $\Vb_n^b(M_1,...,M_n)$ in a normed space $(D,\|\cdot\|_D)$, assume that the map $(M_1,...,M_n) \to \Vb_n^b(M_1,...,M_n)$ is measurable for every $n\in \N$ outer almost surely [the randomness in $\Vb_n^b$ is allowed to come from sources apart from the $M_i$]. Assume that for $i \in \N$ there exist approximations $A_{i,n}^b,A_i$ such that $(M_1,...,M_n) \to A_{i,n}^b(M_1,...,M_n)$ is measurable for every $i,n\in \N$ outer almost surely.
\begin{eqnarray*}
& (\rm{i}) & ~ \text{For every } i\in\N: ~ A_{i,n}^b \weakPi{M} A_{i} ~\text{ for } {n\rightarrow\infty}, \\
& (\rm{ii}) & ~ A_{i} \weak \Vb ~ \text{ for } {i\rightarrow\infty}, \\
& (\rm{iii}) & ~ \text{For every }\varepsilon>0: ~ \lim_{i\rightarrow\infty}\limsup_{n\rightarrow\infty} \Pb^*(\|A_{i,n}^b-\Vb_n^b\|>\varepsilon) = 0.
\end{eqnarray*}
where $A_i,\Vb$ denote a tight processes. Then $\Vb_n^b \weakPi{M} \Vb$.
\end{lemma}

\textbf{Proof of Lemma \ref{lem:bootappr}}
We need to show that
\begin{enumerate}[(a)]
\item  $\sup_{f\in \BL_1} \left| \E_M f(\Vb_{n}^b) - \E f(\Vb) \right| \rightarrow 0$ in outer probability,
\item $\E_M f(\Vb_n^b)^* - \E_M f(\Vb_n^b)_* \Pconv 0$ for all $f\in \BL_1$.
\end{enumerate}
Begin by observing that for every $i \in \N$, every $\omega$, and every $f \in \BL_1$
\begin{eqnarray*}
\left| \E_M f(\Vb_n^b) - \E f(\Vb) \right|
\leq \left| \E_M f(\Vb_n^b) - \E_M f(A_{i,n}^b) \right| + \left| \E_M f(A_{i,n}^b) - \E f(A_i) \right| + \left| \E f(A_i) - \E f(\Vb) \right|.
\end{eqnarray*}
Moreover, for every $\omega$
\begin{eqnarray*}
\sup_{f\in \BL_1} \left| \E_M f(\Vb_n^b) - \E_M f(A_{i,n}^b) \right|
\leq \sup_{f\in \BL_1} \E_M \left|f(\Vb_n^b) - f(A_{i,n}^b)\right|
\leq \E_M \left[\|\Vb_n^b-A_{i,n}^b\|^*\wedge 2\right].
\end{eqnarray*}
In particular, this implies that for any $\gamma>0$
\[
\E^*\Big[ \sup_{f\in \BL_1} \left| \E_M f(\Vb_n^b) - \E_M f(A_{i,n}^b) \right|  \Big]
\leq \E\left[\|\Vb_n^b-A_{i,n}^b\|^*\wedge 2\right]
\leq 2\Pb(\|\Vb_n^b-A_{i,n}^b\|^* > \gamma) + \gamma.
\]
Thus (iii) yields
\[
\lim_{i\rightarrow\infty}\limsup_{n\rightarrow\infty} \E^*\Big[ \sup_{f\in \BL_1} \left| \E_M f(\Vb_n^b) - \E_M f(A_{i,n}^b) \right|  \Big] =  0.
\]
Fix arbitrary $\eps,\eta > 0$. The computations above yield the existence of an $i_1\in\N$ such that for all $i \geq i_1$
\[
\limsup_{n \to \infty}\Pb^* \Big( \sup_{f\in \BL_1} \left| \E_M f(\Vb_n^b) - \E_M f(A_{i,n}^b) \right| > \eps/3 \Big) < \eta/3.
\]
Moreover, by (ii) and the definition of weak convergence, there exists an $i_2\in\N$ such that for all $i \geq i_2$
\[
\Pb^* \Big( \sup_{f\in \BL_1} \left| \E f(A_i) - \E f(\Vb) \right| > \eps/3 \Big) < \eta/3.
\]
Set $k = i_1\vee i_2$. Then (i) implies as $n \to \infty$
\[
\sup_{f\in \BL_1} \left| \E_M f(A_{k,n}^b) - \E f(A_k) \right| = o_{\Pb^*}(1),
\]
and combining all the results above we see that
\[
\limsup_{n\to\infty}\Pb^*\Big( \sup_{f\in \BL_1} \left| \E_M f(\Vb_{n}^b) - \E f(\Vb) \right| > \eps\Big)<\eta.
\]
Since $\eta,\eps$ were arbitrary, this establishes (a). For a proof of (b), note that (i) implies $A_{i,n}^b \weak A_i$ since conditional weak convergence implies unconditional weak convergence [see the proof of Theorem 10.4, assertion $(ii)\Rightarrow (i)$ in \cite{kosorok2008}]. Thus, by an the results in \cite{budevo2011}, (i)-(iii) imply that $\Vb_n^b \weak \Vb$. In particular, this implies asymptotic measurability of $\Vb_n^b$ [see section 1.3 in \cite{vandervaart1996}], and together with the continuity of $f \in \BL_1$ this shows that $f(\Vb_n^b) \weak f(\Vb)$ by an application of the continuous mapping theorem. Thus $\E_M f(\Vb_n^b)^* - \E_M f(\Vb_n^b)_*$ converges to zero in $L^1$, hence also in probability. Now the proof is complete. \hfill $\Box$

\end{appendix}

\bibliography{extreme}

\begin{thebibliography}{}

\bibitem[Adams and Nobel, 2010]{adno2010}
Adams, T. and Nobel, A. (2010).
\newblock Uniform convergence of vapnik--chervonenkis classes under ergodic
  sampling.
\newblock {\em The Annals of Probability}, 38(4):1345--1367.

\bibitem[Andrews and Pollard, 1994]{anpo1994}
Andrews, D.~W. and Pollard, D. (1994).
\newblock {An introduction to functional central limit theorems for dependent
  stochastic processes.}
\newblock {\em International Statistical Review}, 62(1):119--132.

\bibitem[Aue and Reimherr, 2009]{AHR2009}
Aue, A.~Horv{\'a}th, L. and Reimherr, M. (2009).
\newblock Delay times of sequential procedures for multiple time series
  regression models.
\newblock {\em Journal of Econometrics}, 149:174--190.

\bibitem[Berkes et~al., 2009]{beschho2009}
Berkes, I., H{\"o}rmann, S., and Schauer, J. (2009).
\newblock Asymptotic results for the empirical process of stationary sequences.
\newblock {\em Stochastic Processes and Their Applications}, 119(4):1298--1324.

\bibitem[Beutner et~al., 2012]{bewuza2012}
Beutner, E., Wu, W.~B., and Z{\"a}hle, H. (2012).
\newblock Asymptotics for statistical functionals of long-memory sequences.
\newblock {\em Stochastic Processes and their Applications}, 122(3):910--929.

\bibitem[Beutner and Z{\"a}hle, 2010]{beza2010}
Beutner, E. and Z{\"a}hle, H. (2010).
\newblock A modified functional delta method and its application to the
  estimation of risk functionals.
\newblock {\em Journal of Multivariate Analysis}, 101(10):2452--2463.

\bibitem[Beutner and Z{\"a}hle, 2012]{beza2012}
Beutner, E. and Z{\"a}hle, H. (2012).
\newblock Deriving the asymptotic distribution of u-and v-statistics of
  dependent data using weighted empirical processes.
\newblock {\em Bernoulli}, 18(3):803--822.

\bibitem[Billingsley, 1968]{billingsley1968}
Billingsley, P. (1968).
\newblock {\em Convergence of Probability Measures}.
\newblock John Wiley \& Sons Inc., New York.

\bibitem[B{\"u}cher et~al., 2011]{budevo2011}
B{\"u}cher, A., Dette, H., and Volgushev, S. (2011).
\newblock New estimators of the pickands dependence function and a test for
  extreme-value dependence.
\newblock {\em Annals of Statistics}, 39(4):1963--2006.

\bibitem[B{\"u}cher and Ruppert, 2013]{buru2013}
B{\"u}cher, A. and Ruppert, M. (2013).
\newblock Consistent testing for a constant copula under strong mixing based on
  the tapered block multiplier technique.
\newblock {\em Journal of Multivariate Analysis}, 116(0):208 -- 229.

\bibitem[B{\"u}cher and Volgushev, 2011]{buvo2011}
B{\"u}cher, A. and Volgushev, S. (2011).
\newblock Empirical and sequential empirical copula processes under serial
  dependence.
\newblock {\em arXiv preprint arXiv:1111.2778}.

\bibitem[Chu and White, 1995]{CSW1995}
Chu, C-S.~Stinchcombe, M. and White, H. (1995).
\newblock Monitoring structural change.
\newblock {\em Econometrica}, 64(5):1045--1065.

\bibitem[Cs{\"o}rg{\"o} and Horv{\'a}th, 1997]{CH1997}
Cs{\"o}rg{\"o}, M. and Horv{\'a}th, L. (1997).
\newblock {\em Limiting Theorems in Change-point Analysis}.
\newblock New York, Wiley.

\bibitem[Dehling et~al., 2002]{demiso2002}
Dehling, H., Mikosch, T., and S{\o}rensen, M. (2002).
\newblock {\em Empirical process techniques for dependent data}.
\newblock Birkh{\"a}user.

\bibitem[Dehling and Taqqu, 1989]{deta1989}
Dehling, H. and Taqqu, M.~S. (1989).
\newblock The empirical process of some long-range dependent sequences with an
  application to u-statistics.
\newblock {\em The Annals of Statistics}, pages 1767--1783.

\bibitem[Doss and Gill, 1992]{dogi1992}
Doss, H. and Gill, R.~D. (1992).
\newblock An elementary approach to weak convergence for quantile processes,
  with applications to censored survival data.
\newblock {\em Journal of the American Statistical Association},
  87(419):869--877.

\bibitem[Fermanian et~al., 2004]{ferradweg2004}
Fermanian, J.~D., Radulovi{\'{c}}, D., and Wegkamp, M.~H. (2004).
\newblock Weak convergence of empirical copula processes.
\newblock {\em Bernoulli}, 10:847--860.

\bibitem[Gao and Zhao, 2011]{gazh2011}
Gao, F. and Zhao, X. (2011).
\newblock Delta method in large deviations and moderate deviations for
  estimators.
\newblock {\em The Annals of Statistics}, 39(2):1211--1240.

\bibitem[Gill and Johansen, 1990]{gijo1990}
Gill, R.~D. and Johansen, S. (1990).
\newblock A survey of product-integration with a view toward application in
  survival analysis.
\newblock {\em The Annals of Statistics}, pages 1501--1555.

\bibitem[Hagemann, 2012]{hagemann2012}
Hagemann, A. (2012).
\newblock Stochastic equicontinuity in nonlinear time series models.
\newblock {\em Arxiv preprint arXiv:1206.2385}.

\bibitem[Inoue, 2001]{inoue2001}
Inoue, A. (2001).
\newblock Testing for distributional change in time series.
\newblock {\em Econometric Theory}, 17(1):156--187.

\bibitem[Kaplan and Meier, 1958]{kame1958}
Kaplan, E.~L. and Meier, P. (1958).
\newblock Nonparametric estimation from incomplete observations.
\newblock {\em Journal of the American statistical association},
  53(282):457--481.

\bibitem[Kiefer and Vogelsang, 2005]{KV2005}
Kiefer, N.~M. and Vogelsang, T.~J. (2005).
\newblock A new asymptotic theory for heteroskedasticity-autocorrelation robust
  tests.
\newblock {\em Econometric Theory}, 21:1130--1164.

\bibitem[Kitamura, 1997]{Kita1997}
Kitamura, Y. (1997).
\newblock A new asymptotic theory for heteroskedasticity-autocorrelation robust
  tests.
\newblock {\em Annals of Statistics}, 25:2084--2102.

\bibitem[Kosorok, 2008]{kosorok2008}
Kosorok, M.~R. (2008).
\newblock {\em Introduction to Empirical Processes and Semiparametric
  Inference}.
\newblock Springer Series in Statistics, New York.

\bibitem[K\"unsch, 1989]{Kuns1989}
K\"unsch, H. (1989).
\newblock The jackknife and the bootstrap for general stationary observations.
\newblock {\em Annals of Statistics}, 17:1217--1241.

\bibitem[Maejima and Tudor, 2007]{matu2007}
Maejima, M. and Tudor, C.~A. (2007).
\newblock Wiener integrals with respect to the hermite process and a
  non-central limit theorem.
\newblock {\em Stochastic Analysis and Applications}, 25(5):1043--1056.

\bibitem[Perron, 2006]{P2006}
Perron, P. (2006).
\newblock Dealing with structural breaks.
\newblock {\em Palgrave Handbook of Econometrics, vol I: Econometric Theory,
  eds. K. Patterson and T.C. Mills}, pages 278--352.

\bibitem[Politis and Romano, 1994]{PR1994}
Politis, D.~N. and Romano, J.~P. (1994).
\newblock Large sample confidence regions based on subsamples under minimal
  assumptions.
\newblock {\em Annals of Statistics}, 22:2031--2050.

\bibitem[Politis et~al., 1999]{PRW1999}
Politis, D.~N., Romano, J.~P., and Wolf, M. (1999).
\newblock {\em Subsampling}.
\newblock New York, Springer.

\bibitem[Radulovi{\'c}, 2009]{radulovic2009}
Radulovi{\'c}, D. (2009).
\newblock Another look at the disjoint blocks bootstrap.
\newblock {\em Test}, 18(1):195--212.

\bibitem[R{\"u}schendorf, 1974]{ruschendorf1974}
R{\"u}schendorf, L. (1974).
\newblock On the empirical process of multivariate, dependent random variables.
\newblock {\em Journal of Multivariate Analysis}, 4(4):469--478.

\bibitem[R{\"{u}}schendorf, 1976]{rueschendorf1976}
R{\"{u}}schendorf, L. (1976).
\newblock Asymptotic distributions of multivariate rank order statistics.
\newblock {\em Annals of Statistics}, 4:912--923.

\bibitem[Sen, 1974]{sen1974}
Sen, P.~K. (1974).
\newblock Weak convergence of multidimensional empirical processes for
  stationary $\phi $-mixing processes.
\newblock {\em The Annals of Probability}, 2(1):147--154.

\bibitem[Shao, 2010a]{shao2010}
Shao, X. (2010a).
\newblock A self-normalized approach to confidence interval construction in
  time series.
\newblock {\em Journal of the Royal Statistical Society: Series B (Statistical
  Methodology)}, 72(3):343--366.

\bibitem[Shao, 2010b]{shao2010c}
Shao, X. (2010b).
\newblock {\mockalph{zz}}corrigendum: A self-normalized approach to confidence
  interval construction in time series.
\newblock {\em Journal of the Royal Statistical Society: Series B (Statistical
  Methodology)}, 72(5):695--696.

\bibitem[Shao, 2011]{shao2011}
Shao, X. (2011).
\newblock A simple test of changes in mean in the possible presence of
  long-range dependence.
\newblock {\em Journal of Time Series Analysis}, 32(6):598--606.

\bibitem[Shao, 2012]{shao2012}
Shao, X. (2012).
\newblock Parametric inference in stationary time series models with dependent
  errors.
\newblock {\em Scandinavian Journal of Statistics}, 39(4):772--783.

\bibitem[Shao and Politis, 2013]{shpo2013}
Shao, X. and Politis, D. (2013).
\newblock Fixed b subsampling and the block bootstrap: improved confidence sets
  based on p-value calibration.
\newblock {\em Journal of the Royal Statistical Society: Series B (Statistical
  Methodology)}, 75:161--184.

\bibitem[Shao and Zhang, 2010]{shzh2010}
Shao, X. and Zhang, X. (2010).
\newblock Testing for change points in time series.
\newblock {\em Journal of the American Statistical Association},
  105(491):1228--1240.

\bibitem[Van~der Vaart and Wellner, 1996]{vandervaart1996}
Van~der Vaart, A.~W. and Wellner, J.~A. (1996).
\newblock {\em Weak Convergence and Empirical Processes}.
\newblock Springer Verlag, New York.

\bibitem[Vogelsang, 1999]{Vo1999}
Vogelsang, T.~J. (1999).
\newblock Sources of nonmonotonic power when testing for a shift in mean of a
  dynamic time series.
\newblock {\em Journal of Econometrics}, 88:283--299.

\bibitem[Wu, 2007]{Wu2007}
Wu, W. (2007).
\newblock Strong invariance principles for dependent random variables.
\newblock {\em The Annals of Probability}, 35(6):2294--2320.

\bibitem[Wu and Shao, 2004]{wush2004}
Wu, W. and Shao, X. (2004).
\newblock Limit theorems for iterated random functions.
\newblock {\em Journal of Applied Probability}, 41(2):425--436.

\bibitem[Yoshihara, 1975]{yoshihara1975}
Yoshihara, K.-i. (1975).
\newblock Weak convergence of multidimensional empirical processes for strong
  mixing sequences of stochastic vectors.
\newblock {\em Probability Theory and Related Fields}, 33(2):133--137.

\bibitem[Zhou and Shao, 2013]{zhsh2013}
Zhou, Z. and Shao, X. (2013).
\newblock Inference for linear models with dependent errors.
\newblock {\em Journal of the Royal Statistical Society: Series B (Statistical
  Methodology)}, 75:323--343.

\end{thebibliography}

\end{document}